\documentstyle[a4,11pt,amssymb]{article}

\def\section*#1{\begin{center} {\normalsize {\bf #1}} \end{center}}

\font\eulersm=eusm10 at 11pt
\def\esm#1{\hbox{\eulersm {#1}}}

\newcommand{\ee}{\hfill $\Box$}
\newcommand{\Proof}{\noindent {\sl Proof.} \hspace{3mm}}
\newcommand{\vv}{\vspace{2mm}}

\newcommand{\sd}{\rtimes}
\newcommand{\A}{\mbox{\rm Aut}}

\newcommand{\T}{{\esm{T}}}
\newcommand{\ram}{*}

\newcommand{\Ends}{\mbox{{\rm Ends}}}
\newcommand{\B}{\Sigma}
\newcommand{\an}{\mbox{{\rm {\tiny an}}}}
\begin{document}

{\footnotesize \noindent Running head: Automorphisms of Mumford curves \hfill August 15, 2000

\noindent Math.\ Subj.\ Class.\ (1991): 14E09, 20E08, 30G06  \hfill {\bf Accepted version}\\
\noindent \mbox{} \hfill {\bf Correction added Sept 17. 2019}}

\vspace{3cm}

\begin{center} 

{\Large Discontinuous groups in positive characteristic

\vv

and automorphisms of Mumford curves}

\vv

{\sl by} Gunther Cornelissen, Fumiharu Kato {\sl and} Aristides Kontogeorgis

\vv

{\footnotesize published Math.\ Ann.\ {\bf 320}, 55--85 (2001);\\ correction published https://doi.org/10.1007/s00208-019-01901-9 (2019)}

\end{center}

\vv

\vv

\section*{Introduction}

\vv

\noindent All compact Riemann surfaces ({\sl vis-\`a-vis} complex projective curves) of genus $g \geq 2$ share the same universal topological covering space, and hence admit a uniformization by a discrete subgroup of $PSL(2,{\bf R})$. This description  leads to a fruitful dualism in the theory of complex curves, visible in both its theorems and their proofs. The most appropriate instance of this to be quoted here is probably the fact that Hurwitz's upper bound $84(g-1)$ on the number of automorphisms of a compact Riemann surface  of genus $g \geq 2$ is equivalent to Siegel's lower bound $\pi / 21$ on the volume of the fundamental domain of a Fuchsian group ({\sl cf.} Lehner, \cite{Lehner:64}, p. 402-405). 

Although over a {\sl non-archimedean} valued field $k$, compact rigid-analytic and projective algebraic curves can still be identified, the analogue of the dualism to discrete group theory is of quite a different nature, as not all curves admit the same universal topological covering space ({\sl e.g.}, curves having good reduction are analytically simply connected). But Mumford has shown that curves whose stable reduction is split multiplicative ({\sl i.e.,} a union of rational curves intersecting in $k$-rational points) are isomorphic to an analytic space of the form $\Gamma \backslash ({\bf P}^1_k-{\cal L}_\Gamma)$, where $\Gamma$ is a discontinuous group in $PGL(2,k)$ with ${\cal L}_\Gamma$ as set of limit points. Thus, the theory of non-archimedean discrete groups is both more restrictive than the complex analytic one (as it cannot be applied to any curve), and more powerful, as it can lead to stronger results for such so-called Mumford curves. 

As an example of this, Herrlich has shown that for $p$-adic Mumford curves of genus $g \geq 2$, Hurwitz's bound can be strengthened to  $12(g-1)$ (if $p \geq 7$, {\sl cf.} \cite{Herrlich:80}). One can wonder what happens for Mumford curves
over non-archimedean fields of {\sl positive characteristic}. 

It ought to be mentioned that in positive characteristic, general algebraic curves can have many more automorphisms than expected from Hurwitz's bound, {\sl e.g.}, there exist algebraic curves of arbitrary high genus with more than $16g^4$ automorphisms (Stichtenoth, \cite{Stichtenoth:73}). It has been observed on many occasions that the most anomalous examples invariably have a low $p$-rank. One is therefore led to expect better properties for so-called ordinary curves, {\sl i.e.,} curves $X$ of genus $g$ for which 
$$ \dim_{{\bf Z}/p} \mbox{Jac}(X)[p] = g, $$
which are Zariski-dense in the moduli space of curves of genus $g$.
Indeed, S.\ Nakajima (\cite{Nakajima:87}) has shown that for such curves, $|\mbox{Aut}(X)| \leq 84g(g-1)$. However, one knows of no infinite set of ordinary curves of arbitrary high genus all of whose elements attain this bound. Rather, the number of automorphisms in all such known collections does not exceed a cubic polynomial in $\sqrt{g}$, and it has been suggested that worse cannot happen. 

Mumford curves are known to be ordinary (\S 1), whence the aforementioned pathologies are moderated by Nakajima's result. It is our aim to show that a bound of the type suggested before indeed holds for Mumford curves:

\vv

{\bf Theorem.} \ \ {\sl Let $X$ be a Mumford curve of genus $g \geq 2$ over a non-archimedean valued field of characteristic $p>0$. Then
$$|\A(X)| \leq \max \{ 12(g-1), 2 \sqrt{g} (\sqrt{g}+1)^2\}. $$}

Actually, the proof yields more: it provides a kind of classification of those curves $X$ for which $|\A(X)|>12(g-1)$ (\S 4). 

\vv

{\bf Remark.} \ \ A reformulation of the above theorem in the style of Siegel's lower bound is as follows: {\sl if \ $\Gamma$ is a Schottky group of rank $g$, then the $\mu$-invariant of its normalizer $N$ {\rm (cf.\ \cite{Karrass:73}, \cite{Herrlich:80})} is bounded from below by 
$$ \mu(N) \geq \min \{ \frac{1}{12}, \frac{\sqrt{g}-1}{2\sqrt{g}(\sqrt{g}+1)} \}. $$}

The proof of the main theorem uses the description of the automorphism group by
the normalizer $N$ of its Schottky group $\Gamma$, whose structure can be studied by 
its action on the Bruhat-Tits tree $\T$ of $PGL(2,k)$. Although this principle of proof was also
used by Herrlich in the $p$-adic case, the quintessence of our techniques is
rather different. In the $p$-adic case the expected bound is linear in the
genus, and this allows Herrlich to restrict to normalizers that are the 
amalgam of {\sl two} finite groups. In positive characteristic, the expected bound is not linear in the genus, so we have to consider more
complicated normalizers. We therefore investigate directly the link between the ramification in $\pi\ : \ X \rightarrow S:= \A(X)\backslash X$ and the combinatorial geometry of the analytic reduction of $S$. From Hurwitz's formula applied to $\pi$, it is immediate that the bound $12(g-1)$ on the order of the automorphism group holds, unless  $S={\bf P}^1$ and $\pi$ is branched above $m=2$ or $m=3$ points. Let $\T_N$ (respectively ${\T}_N^\ram$) be the subtree of $\T$ which is generated by the limit points of $\Gamma$, seen as ends of $\T$ (respectively the limit points together with the fixed points of torsion elements in $N$). The quotients $T_N=N \backslash {\T}_N$ and $T_N^\ram = N \backslash {\T}_N^\ram$ are homotopic to the intersection graph of the analytic reduction of $S$ as a rigid analytic space, but $T_N^\ram$  has finitely many ends attached to it, which are in bijection with the points above which $\pi$ is branched. The advantage of using $T_N^\ram$ instead of the usual $T_N$ lies in the following key proposition, which replaces the ``restriction to 
an amalgam of two groups'' in the $p$-adic case: 

\vv

{\bf Proposition 1.} \ \ {\sl If \ $T$ is a subtree of $T_N^\ram$ having the same ends as $T_N^\ram$, then $T$ is a contraction of $T_N^\ram$, {\rm i.e.,} every geodesic connecting a point from $T_N^\ram-T$ to $T$ is a path on which the stabilizers of vertices are ordered increasingly w.r.t.\ inclusion in the direction of $T$. In particular, the amalgams associated to $T$ and $T_N^\ram$ are the same $(=N)$.}

\vv

The proof of this result is very combinatorial and depends on the structure 
theorem for finite subgroups of $PGL(2,k)$ (\S 3).

By applying proposition 1, it is enough to consider a simpler subtree $T$ of $\T_N^\ram$ which is a ``line'' if $m=2$ and a ``star'' if $m=3$. We then show that only finitely many types of such trees $T$ (and hence, of such groups $N$) exist (\S 4). For each such type of group, we use the link between the ends of
the corresponding $T$ and the ramification in $\pi$ to prove the desired bound (\S 5).  

\vv 

{\bf Remark.} \ \ Proposition 1 holds in any characteristic, and can be applied to obtain a classification of Galois covers of ${\bf P}^1_{{\bf C}_p}$ which are Mumford curves ramified above three points, as in \cite{Kato:99}. 

\vv

We now turn to some comments on the main theorem. Note that $12(g-1)> 2 \sqrt{g} ( \sqrt{g}+1)^2$ only for $g \in \{5,6,7,8 \}$. But there do exist curves of such genus for which the bound is attained. Let $M_g$ (respectively, ${\cal M}_g$) denote the moduli space of algebraic curves of genus $g$ (respectively, of Mumford curves of genus $g$), considered as algebraic space (respectively, rigid analytic space). We then have the following result, which identifies a very special locus in $M_6$:

\vv

{\bf Proposition 2.} \ \ {\sl If $X$ is a Mumford curve of genus $g$ over a field $k$ of characteristic $p>0$ such that $\A(X)$ is a non-solvable group of order $>2 \sqrt{g} (\sqrt{g}+1)^2$, then $p \neq 2,5$ and $X$ is a curve of genus 6 with an icosahedral isomorphism group. For $p>5$, there is a  one-dimensional open stratum $M_6$, along which the automorphism group is icosahedral. It intersects ${\cal M}_6$ in a
one-dimensional rigid analytic open set} (\S 7). 

\vv

Let $\ast$ denote amalgamated product. The bound in the above theorem is sharp in the following sense: 

\vv

{\bf Proposition 3.} \ \ {\sl For fixed $p$, there exist Mumford curves $X_{t,c}$ of arbitrary high genus $g_t = (p^t-1)^2$ that attain the bound in the main theorem; $X_{t,c}$ is the projective curve with affine equation
$ (x^{p^t}-x)(y^{p^t}-y) = c, $
for some $c \in k^*$ with $|c|<1$. The normalizer of its
Schottky group is isomorphic to 
$$N_t = ({\bf Z}_p^t \sd {\bf Z}_{p^t-1}) \ast_{{\bf Z}_{p^t-1}} D_{p^t-1},$$ where 
${\bf Z}_n,D_n$ are the cyclic and dihedral group of order $n$ and $2n$ respectively. 
The Schottky group of $X_{t,c}$ is generated by the commutators $[\varepsilon, \gamma \varepsilon \gamma]$, for all $\varepsilon \in   
{\bf Z}_p^t$, where $\gamma$ is a fixed involution in $D_{p^t-1}-{\bf Z}_{p^t-1}$. For varying $c$, the $X_{t,c}$ form a one-dimensional open rigid analytic stratum of curves in ${\cal M}_{g_t}$ that attain the bound} (\S 9).

\vv

Another less obvious family of Mumford curves whose number of automorphisms exceeds Hurwitz's bound is given by moduli schemes for rank two Drinfeld modules with principal level structure ({\sl e.g.,} Gekeler \& Reversat \cite{Gekeler:96}). The analytic description is as follows: Let $q=p^t$, $F={\bf F}_q(T)$, and $A={\bf F}_q[T]$; let 
$F_\infty={\bf F}_q((T^{-1}))$  be the completion of $F$ and $C$ a completion of the algebraic closure of $F_\infty$. On Drinfeld's ``upper half plane'' $\Omega:= {\bf P}^1_C - {\bf P}^1_{F_\infty}$ (which is a rigid analytic space over $C$), the group $GL(2,A)$ acts by fractional transformations. Let ${\cal Z} \cong {\bf F}_q^*$ be its center. For ${\frak n} \in A$, the quotients of $\Omega$ by congruence subgroups $\Gamma({\frak n}) = \{ \gamma \in GL(2,A) : \gamma = {\bf 1} \mbox{ mod } {\frak n} \}$ are open analytic curves which can be
compactified to projective curves $X({\frak n})$ by adding finitely many cusps. They turn out to be Mumford curves for the free group which splits the inclusion $\Gamma({\frak n})_{\mbox{\tiny tor}} \lhd \Gamma({\frak n})$, where $\Gamma({\frak n})_{\mbox{\tiny tor}}$ is the subgroup generated by torsion elements ({\sl cf.} Reversat \cite{Reversat:96}).

\vv

{\bf Proposition 4.} \ \ {\sl Let $p \neq 2$ or $q \neq 3$. The automorphism group of the Drinfeld modular 
curves $X({\frak n})$ is the ``modular'' automorphism group $G({\frak n}):=\Gamma(1)/\Gamma({\frak n}) {\cal Z}$. The normalizer of its Schottky group is isomorphic to 
$$ PGL(2,p^t) \ast_{{\bf Z}_p^t \sd {\bf Z}_{p^t-1}} {\bf Z}_p^{td} \sd {\bf Z}_{p^{t}-1}, $$
which only depends on $d=\deg({\frak n})$. For fixed $g$ with $g=g(X({\frak n}))$ for some $n$, $X({\frak n})$ is one of at most finitely many Mumford curves with $G({\frak n})$ as automorphism group} (\S 10).

\vv

Note that, whereas the corresponding statement for classical modular curves is easy to prove using the Hurwitz formula (\cite{Serre:96}), the occurrence of wild ramification in $X({\frak n}) \rightarrow X(1)$ prevents us from giving an easy proof of the above result. Note also that the final statement of proposition 4 highlights the exceptional status (as a point of the moduli space) of such modular curves which admit {\sl simultaneous} uniformization by a Schottky group and by an arithmetic group on an open part. 

We want to finish this introduction by mentioning a few questions raised by 
this work. 

\vv

{\bf Question 1} (ordinary curves). \ \ Does the above bound with $12(g-1)$ replaced by $84(g-1)$ hold for all
ordinary curves? 

\vv

{\bf Questions 2} (Drinfeld modular curves). \ \ Can one be more precise about the abstract structure of the Schottky group of $X({\frak n})$? It should depend, not only on $\deg({\frak n})$, but also on the decomposition
of ${\frak n}$ into prime ideals of $A$.

\vv

\vv

\section*{1. Automorphisms and uniformization in positive characteristic}

\vv 

(1.1) Let $X$ be a complete irreducible curve of genus $g \geq 2$ over a field $k$ of characteristic $p > 0$. Let $K$ denote the algebraic closure of $k$, and $\bar{k}$ its residue field. Let
$\gamma$ denote the $p$-rank of $X$, {\sl i.e.,} $\gamma=\dim_{{\bf Z}/p} \mbox{Jac}(X)(K)[p]$. If $\gamma=g$, then $X$ is called {\sl ordinary} (such curves are dense in the moduli space of curves of genus $g$). Recall (\cite{Gerritzen:80}, \cite{Mumford:72}) that $X$ is called a {\sl Mumford
curve} if it is ``uniformized over $k$ by a Schottky group''. This means that (a) $k$ is complete with respect to a non-archimedean valuation; (b) there exists
a free group $\Gamma$ of rank $g$ in $PGL(2,k)$, acting on ${\bf P}^1_k$ with limit set ${\cal L}_\Gamma$; (c) the analytification $X^{\an}$ of $X$ satisfies $X^{\an} \cong \Gamma
\backslash ({\bf P}^1_k - {\cal L}_\Gamma)$ as rigid analytic spaces. Mumford has shown that these conditions are equivalent to the stable reduction of $X$ only having rational components with $\bar{k}$-rational double points. Because of the ``GAGA''-correspondence for one-dimensional rigid analytic spaces, we do not have to (and will not) distinguish between analytic and algebraic curves. The following result is folklore:

\vv

{\bf (1.2) Lemma.} \ \ {\sl Mumford curves are ordinary.}

\vv

\Proof It is known by the works of Manin-Drinfeld (\cite{Manin:73}) and Gerritzen (\cite{Gerritzen:72}) that the Jacobian of a Mumford curve $X$ associated to a free group
$\Gamma$ can be uniformized as a rigid
analytic abelian variety by a ``multiplicative lattice'', {\sl i.e.,} 
$$ \mbox{Jac}(X)(K) \cong (K^*)^g/\Lambda $$
for some lattice $\Lambda$ (actually, $\Lambda \cong \Gamma^{\mbox{{\tiny ab}}}$) with multiplicative basis $\{ \lambda_1,...,\lambda_g \}$. Then $$ \mbox{Jac}(X)[p] \cong \sqrt[p]{\Lambda}/\Lambda, $$
where $\sqrt[p]{\Lambda}$ is generated by $\{ \sqrt[p]{\lambda_1},...,\sqrt[p]{\lambda_g} \}$ (since the $p$-th roots of unity are trivial in ${k}$), which clearly has ${\bf Z}/p$-rank $= g$. 
\ee

\vv

The following is the rigid analytic analogue of a well known
theorem on conformal automorphisms of Riemann surfaces ({\sl e.g.,} \cite{Lehner:64}, VI.3.L):

\vv

{\bf (1.3) Theorem} (Gerritzen \cite{Gerritzen:80}, VII.1). \ \ {\sl Let $X$ be a Mumford curve 
over $k$ uniformized by a Schottky group $\Gamma$. Let $N$ be the normalizer of $\Gamma$ in 
$PGL(2,k)$. Then $\A(X) \cong N/\Gamma$. \ee}

\vv

Although $N$ is not free, it makes sense to consider the 1-dimensional
analytic quotient space $S:= N\backslash ({\bf P}^1_k-{\cal L}_N)$, which can be algebraized, since it is the quotient of $X$ by the finite group $\A(X)=N/\Gamma$.

\vv

{\bf (1.4) Notation.} \ \  Denote the genus of $S$ by $r$, and denote by $m$ the number of points of $S$ above which $\pi : X \rightarrow S$ is branched. 

\vv

{\bf (1.5) Lemma.} \ \ {\sl If $r>0$ or $m>3$, then $|\A(X)| \leq 12(g-1)$. Furthermore, the case $r=0,m=1$ cannot occur.} 

\vv

\Proof The first part is an easy consequence of the Riemann-Hurwitz formula applied to $\pi$.
({\sl cf.} 
Stichtenoth's proof of {\sl Satz} 3  in \cite{Stichtenoth:73}). The second part
follows from (1.2) and the same formula, since S.\ Nakajima has shown in \cite{Nakajima:87}, Theorem 2, that for ordinary $X$, the ramification groups of order $\geq 2$ in $\pi$ vanish. \ee

\vv

\vv

\section*{2. Structure of $T_N^\ram$}

\vv 

From now on, we will assume that $r=0$ and give a description of the map $\pi$ by looking at its analytic reduction.

\vv

(2.1) {\sl The Bruhat-Tits tree.} \ \
By taking a finite extension of $k$ if necessary, we can assume that all fixed points of $N$ acting on ${\bf P}^1_k$ are $k$-valued ({\sl cf.} \cite{Gerritzen:80}, I.3.3.2).
Let
$\T$ denote the Bruhat-Tits tree of $PGL(2,k)$ ({\sl i.e.,} its vertices are similarity classes $\Lambda$ of rank two $\cal O$-lattices in $k^2$, and two vertices are connected  by an edge if the corresponding quotient module has length one -- {\sl see} Serre \cite{Serre:80}, Gerritzen \& van der Put \cite{Gerritzen:80}). It is a regular
tree in which the edges emanating from a given vertex are in one-one correspondence with ${\bf P}^1(\bar{k})$, where $\bar{k}$ is the residue field of $k$. The tree $\T$ admits a left action by $PGL(2,k)$, and hence also by any subgroup $G \subseteq N$; by our assumption on the fixed points, $G$ acts without inversions on $\T$.

\vv

{\bf (2.2) Notations.} For any subtree $T$ of $\T$, let $\Ends(T)$ denote its set of ends ({\sl i.e.,}
equivalence classes of half-lines differing by a finite segment). There is a natural correspondence between ${\bf P}^1(k)$ and  $\Ends(\T)$.  Let $V(T)$ and $E(T)$ denote the set of vertices and edges of $T$ respectively. For any $\Lambda \in V(T), \sigma \in E(T)$, let $\sigma \dashv \Lambda$ mean that $\sigma$  originates at $\Lambda$.  Let $G_x$ denote the stabilizer of a vertex or edge $x$ of $T$ for the action of $G$. Vertical bars $| \cdot |$ will denote
geometric realization. For any $u,v \in {\bf P}^1(k)$, let $]u,v[$ denote
the apartment in $\T$ connecting $u$ and $v$ (seen as ends of $\T$). Furthermore, for a vertex $\Lambda \in V(\T)$, let $[\Lambda,u[$ denote the half-line connecting $\Lambda$ to $u$ (seen as an end of $\T$). The natural metric on $|\T|$ is denoted by $d$. 

\vv

(2.3) {\sl Action of elements of finite order.} \ \ Let $G$ be a subgroup of $N$. If $\gamma \in G$ is an element of finite
order, let the {\sl mirror} $M(\gamma)$ of $\gamma$ be defined as the smallest subtree 
of $\T$ generated by pointwise $\gamma$-fixed vertices of $\T$ (this definition differs from the one given in \cite{Kato:99}).

If $\gamma$
is an elliptic element ({\sl i.e.}, an element in $G$ of finite order having two distinct eigenvalues of the same valuation), then it has two fixed points in ${\bf P}^1(k)$, and $M(\gamma)$ is just the geodesic connecting them (recall that $k$ is of positive characteristic). 

If $\gamma$ is a parabolic element ({\sl i.e.}, an element in $G$ of finite order
having only one eigenvalue), then it has a unique fixed point $z$. Assume that $k$ is discretely valued and $\bar{k}$ is finite. Let $\cal O$ be the ring of integers of $k$ and $\varpi$ a uniformizer of $k$. For every $w \in {\bf P}^1(k)$ different from $z$, there is a unique vertex $\Lambda_w \in ]w,z[$ such that $M(\gamma) \cap ]w,z[= [\Lambda_w,z[$. For any two $w,w'$, 
$\Lambda_w$ and $\Lambda_{w'}$ lie at the same distance from $[\Lambda_w,z[\cap[\Lambda_{w'},z[$. Indeed, let $\{e_0,e_1\}$ be the standard basis for $k^2$, and let the parabolic element $\gamma$ be a translation $\gamma: X \mapsto X+t$ for some $t \in k^*$, so that its fixed point $z$ is $z=\infty$. Let $w=0$. Vertices in the apartment connecting $0$ and $\infty$ are similarity classes of lattices $\Lambda_j = {\cal O} e_0 + {\cal O} \varpi^j e_1$, and a short computation shows that $\Lambda_j \in M(\gamma)$ if and only if $j \geq -v(t)$. 
\unitlength0.7pt
\begin{center}
\begin{picture}(150,70)
\multiput(0,10)(40,0){3}{\circle*{4}}
\multiput(0,10)(40,0){2}{\line(1,0){40}}
\multiput(40,10)(40,0){2}{\line(1,1){30}}
\multiput(70,40)(40,0){2}{\circle*{4}}
\put(110,40){\line(0,1){40}}
\put(110,40){\line(2,1){36}}
\multiput(110,40)(0,40){2}{\circle*{4}}
\put(146,58){\circle*{4}}
\put(80,10){\line(1,0){20}}
\multiput(110,10)(10,0){3}{\circle*{1.4}}
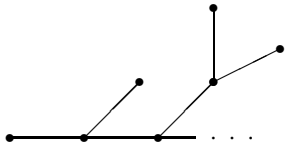
\put(-90,-15){{\footnotesize Picture 1. The mirror of a parabolic element if $q=2$}}
\end{picture}
\end{center}

\vv

(2.4) {\sl Trees associated to subsets of ${\bf P}^1(k)$.} \ \ We can construct a locally finite tree  $\T({\cal L})$ (possibly empty) from any compact subset ${\cal L}$ of ${\bf P}^1(k)$: it is the minimal subtree of $\T$ whose set of ends coincides with ${\cal L}$, or equivalently, the minimal subtree of $\T$ containing  
$ \bigcup_{u,v \in {\cal L}} \ ]u,v[. $

For any subgroup $G$ of $N$ we define $\T_G$  to be the tree associated to the subset ${\cal L}_G$ consisting of the limit points of $G$ in ${\bf P}^1(k)$. If $G$ is finite, $\T_G$ is empty. On the other hand, if $G$ is a finitely generated discrete group, then $\T_G$ coincides with the tree of $G$ as it is defined in Gerritzen \& van der Put \cite{Gerritzen:80}. We also define $\T_G^\ram$ to be the tree associated to the set consisting of all limit points of $G$ {\sl and} all fixed points of elements of finite order in $G$. 

$\T_G$ and $\T_G^\ram$ admit a natural action of $G$, and we denote the quotient graphs by
$T_G:= G \backslash \T_G$ and $T_G^\ram:= G \backslash \T_G^\ram$ respectively, and both corresponding quotient maps will be written as $\pi_G$, by slight abuse of notation. If they are trees, we can consider $T_G$ and $T_G^\ram$ as subtrees of $\T_G$ and $\T_G^\ram$ by choosing any section of $\pi_G$. 

Let $S$ be as in (1.4). The dual graph of the analytic reduction of the 
curve $S$ (with the admissible covering arising from the standard one on 
${\bf P}^1_k - {\cal L}_N$) is exactly the quotient graph $T_N$. Since $S$ is assumed to have genus $r=0$, $T_N$ has no cycles, and hence is a tree. Since $S$ is compact, $T_N$ is finite. 

The tree $T_N^\ram$ has finitely many ends which are in one-to-one correspondence with the branch points in $\pi : X \rightarrow S$. Actually, the stabilizer of such an end ({\sl viz.,} the eventual stabilizer of edges and vertices in the end) is exactly the ramification group of the corresponding branching point in $S$. Note that since
$X$ is ordinary, only the first ramification groups in $\pi : X \rightarrow S$ are non-trivial (Nakajima, {\sl loc.\ cit.}), and hence the stabilizer of an end in $T_N^\ram$ is 
a semi-direct product of a cyclic group and an elementary abelian $p$-group by
the structure theorem for decomposition groups in global field extensions 
({\sl cf.} Serre \cite{Serre:68}, IV \S 2). 

\vv

{\bf (2.5) Example.} \ \  If $G$ is a cyclic group generated by an elliptic element $\gamma$, then
$T_G^\ram=M(\gamma)$. 

\vv

(2.6) {\sl Graphs of groups.} \ \ We turn both of $T_N$ and $T_N^{\ram}$ into graphs of groups by labeling every vertex $v \in V(T_N)$ and edge $e=[vw] \in E(T_N)$ with the stabilizer of its lifts (by any fixed section of $\pi_N$) $\Lambda \in V(\T_N)$ and $\sigma \in E(\T_N)$ for the action of $N$, denoted by $N_v$ and $N_e = N_{vw}=N_v \cap N_w$ respectively. Note that these groups are finite since $N$ is discrete. Then $N$ equals the ``tree product''
of $T_N$, {\sl viz.,} $N$ is the amalgam of the $N_v$ along the $N_e$ (Serre, \cite{Serre:80}, 4.4). In view of what follows in section 4, let us give two examples of how to deal explicitly with such amalgams:

\vv

{\bf (2.6.1) Example.} If $T_N$ is just a segment,{\sl viz.}, 

{\footnotesize
\unitlength0.7pt
\begin{center}
\begin{picture}(40,30)
\multiput(0,15)(35,0){2}{\circle*{4}}
\multiput(0,15)(35,0){1}{\line(1,0){35}}
\put(-45,12){${T_{N}}=$}
\put(-6,21){$G_1$}
\put(30,21){$G_2$}
\put(10,3){$H$}
\put(50,12){,}
\end{picture}
\end{center}}

\noindent then $N=G_1 *_H G_2$ is the ordinary ``product with amalgamation'', {\sl viz.}, if $\langle g_i^\alpha | R_i \rangle$ is a presentation for $G_i$, and $\phi_i : H \rightarrow G_i$ are the occurring injections,  then $N$ admits
the presentation $ N = \langle g_1^\alpha, g_2^\beta \ | \ R_1, R_2,  \phi_1(h)\phi_2(h)^{-1}, \ \forall h \in H \rangle.$ 

From such explicit presentations, it is easy to deduce an isomorphism $$(G_1 \sd H) *_H (G_2 \sd H) = (G_1 * G_2) \sd H$$ for any finite groups $G_1,G_2$ and $H$. Such an isomorphism is used in the proof of (6.8). 

Also observe that, if ${T^{'}_{N'}}$ is the following tree: 

{\footnotesize
\unitlength0.7pt
\begin{center}
\begin{picture}(80,30)
\multiput(0,15)(35,0){3}{\circle*{4}}
\multiput(0,15)(35,0){2}{\line(1,0){35}}
\put(-45,12){${T^{'}_{N'}}=$}
\put(-6,21){$G_1$}
\put(30,21){$H'$}
\put(66,21){$G_2$}
\put(10,3){$H'$}
\put(45,3){$H$}
\put(90,12){,}
\end{picture}
\end{center}}

\noindent then $N=N'$ although $T_N \neq T^{'}_{N'}$. Thus, $T_N$ determines $N$ uniquely, but the converse is not true. This is why, later on, we will classify $N$ instead of $T_N$, and also why we denote the group $N$ by a slightly different diagram
(using amalgamation signs) as the tree (using realizations as plane graphs). \ee

\vv

{\bf (2.6.2) Example.} Suppose that $T_N$ is a ``star'':

{\footnotesize
\unitlength0.7pt
\begin{center}
\begin{picture}(80,70)
\multiput(0,35)(35,0){2}{\circle*{4}}
\put(0,35){\line(1,0){35}}
\put(35,35){\line(1,1){30}}
\put(35,35){\line(1,-1){30}}
\multiput(66,66)(0,-62){2}{\circle*{4}}
\put(-65,30){$T_N= \ \ G_1$}
\put(45,30){$G_2$}
\put(8,23){$H_{12}$}
\put(74,63){$G_3$}
\put(74,1){$G_4$}
\put(30,54){$H_{23}$}
\put(30,10){$H_{24}$}
\put(100,30){.}
\end{picture}
\end{center}}

\noindent We denote the corresponding group $N$ by a similar diagram, in which the edges 
have been replaced by $\ast$-symbols ({\sl cf.\ }(4.7), where such tree products
occur). 
The above picture means that, if $\langle g_i^\alpha | R_i \rangle$ is a presentation for $G_i$, and $\phi^1_{ij} : H_{ij} \rightarrow G_i,\ \phi^2_{ij} : H_{ij} \rightarrow G_j$ are the occurring injections, then $N$ admits a presentation of the form
$$ N = \langle g_1^\alpha, g_2^\beta, g_3^\gamma | R_1, R_2, R_3, \phi^1_{ij}(h_{ij})\phi^2_{ij}(h_{ij})^{-1}, \  \forall h_{ij} \in H_{ij}    \rangle. \Box $$

\vv

The following observation is easy, but it allows us to deduce properties
about the local structure of our tree from information about finite groups
acting on the projective line: 

\vv 

{\bf (2.7) Lemma.} \ \ {\sl Let $v$ be any vertex of $T_N$ corresponding to $\Lambda$ in $\T_N$ and $N_v$ its stabilizer. The canonical bijection $\{ \sigma \in \T : \sigma \dashv  \Lambda \}  \stackrel{\sim}{\longrightarrow} {\bf P}^1(\bar{k})$ induces  a  representation $\rho : N_v \rightarrow PGL(2,\bar{k})$.  Then for every subgroup $G \subseteq N_v$ there is a bijection
between $\{ \sigma \in E(\T^\ram_N) : \sigma \dashv \Lambda \mbox{ {\rm and} } N_\sigma \supseteq G \}$ and points on ${\bf P}^1(\bar{k})$ fixed by $\rho(G)$. Taking quotients by $N$, this leads to a bijection
 $$\left\{   \begin{minipage}{50mm}\setlength{\baselineskip}{.85\baselineskip} {\small  $e \in E(T_N^\ram)$ with  $ e \dashv v$ and  $e$ is fixed by an $N_v$-conjugate of $G$ } \end{minipage} \right \} \stackrel{\sim}{\longrightarrow} \left\{  \begin{minipage}{50mm}\setlength{\baselineskip}{.85\baselineskip} {\small  $ x \in  \rho(N_v) \backslash {\bf P}^1(\bar{k})$ such that points above $x$ in ${\bf P}^1(\bar{k})$ are fixed by a $\rho(N_v)$-conjugate of $\rho(G)$ } \end{minipage} \right\} . \Box$$ }

\vv

We now turn to a description of the possible $N_v$. 

\vv

{\bf (2.8) Notation.}  \ We let ${\bf Z}_n$ denote the cyclic group of
order $n$, $D_n$ the dihedral group of order $2n$ and $E_n$ the elementary abelian $p$-group
${\bf Z}_p^n$ of order $p^n$. Let $T=A_4,O=S_4,I=A_5$ be the tetrahedral, octahedral and icosahedral groups respectively. We also introduce the following
short hand notation: $B(t,n):= E_t \sd {\bf Z}_n$ for $n|p^t-1$; if $n=p^t-1$, this is a Borel subgroup of $PGL(2,p^t)$. Let $\langle \gamma \rangle$ denote the cyclic group generated by $\gamma$. Groups not containing a $p$-group will be called {\sl classical.}

We will write $P(2,q)$ to denote either $PGL(2,q)$ or $PSL(2,q)$ by slight abuse of notation, with the convention that any related numerical 
quantities that appear between set delimiters $\{ \}$ are only to be considered for $PSL(2,q)$.

\vv

We recall the classification of finite subgroups of $PGL(2,k)$, due to L.E.\ Dickson ({\sl cf.} Huppert, \cite{Huppert:79}, II.8.27). A geometric formulation, which is more convenient for our purpose, is due to Valentini and Madan (\cite{Madan:80}). 
From our point of view, the geometric form of the theorem describes the structure of $T_G^\ram$ for finite subgroups $G$ of $PGL(2,k)$.

\vv

{\bf (2.9) Theorem.} \ \ {\sl Any finite subgroup of $PGL(2,k)$ is isomorphic to a finite subgroup of $PGL(2,p^m)$ for some $m>0$. The group $PGL(2,p^m)$ has the following finite subgroups $G$, such that $\pi_G$ is branched over $d$ points with ramification groups isomorphic to
$G_1,...,G_d$:

(i) $G={\bf Z}_n$ for $(n;p)=1$, $d=2$, $G_1=G_2={\bf Z}_n$;

(ii) $G=D_n$ with $p \neq 2$, $n|p^{m} \pm 1$, $d=3$, $G_1=G_2={\bf Z}_2,G_3={\bf Z}_n$ or also, $p=2$, $(n;2)=1$, $d=2$ and $G_1={\bf Z}_2,G_2={\bf Z}_n$; 

(iii) $G=B(t,n)$ for $t \leq m$ and
$n|p^m-1, n|p^t-1$ with $d=2$ and $G_1=G,G_2={\bf Z}_n$ if $n>1$ and $d=1,G_1=G$ otherwise;

(iv) $G=P(2,p^t)$ with $d=2$ and $G_1=B(t,\{\frac{1}{2}\}(p^t-1)),G_2={\bf Z}_{\{\frac{1}{2}\}(p^t+1)}$; 

(v) $T$ of $p \neq 2,3$, $d=3$, $G_1={\bf Z}_2, G_2=G_3={\bf Z}_3$;

(vi) $O$ if $p \neq 2,3$, $d=3$, $G_1={\bf Z}_2, G_2=G_3={\bf Z}_4$;

(vii) $I$ if $5|p^{2m}-1$ and $p \neq 2,3,5$ with $d=3$ and $G_1={\bf Z}_2,G_2={\bf Z}_3,G_3={\bf Z}_5$ or $p=3$, $d=2$ and $G_1=B(1,2),G_2={\bf Z}_5$. \ee}

\vv

{\bf (2.10) Lemma.} \ \ {\sl Let $v \in V(T_N^\ram)$ with preimage $\Lambda \in V(\T_N^\ram)$. The group representation $\rho : N_v \rightarrow PGL(2,k)$ arising from the action of $N_\Lambda$ on its neighbouring vertices in $\T$ is faithful, except possibly when $N_v=B(t,n)$, and if so, ker$(\rho)$ is a $p$-group.}

\vv

\Proof An element $\gamma \in \ker(\rho)$ fixes all edges $\sigma \dashv \Lambda$, hence all neighbouring vertices of $\Lambda$ are in its mirror 
$M(\gamma)$. But if $\gamma$ is elliptic, then its mirror consists of only one apartment, whence cannot stabilize more than two non-collinear vertices in $\T$, whereas the valency of $\T$ is $q+1>2$. So $\gamma$ has to be parabolic, and $\ker(\rho)$ is a normal $p$-subgroup of $N_\Lambda$. Since $PSL(2,q)$ is simple, it follows from (2.9) that the only finite subgroups of $PGL(2,k)$ which have such a non-trivial normal $p$-subgroup are of the form $B(t,n)$. \ee

\vv

\vv

\section*{3. Contraction of $T_N^\ram$ -- proof of proposition 1}

\vv

{\bf (3.1) Definition.} \ \ Let $T$ be a tree of groups in $\T$. A subtree of groups
$T'$ of $T$ is called a {\sl contraction} of $T$ if every geodesic connecting a point from $T-T'$ to $T'$ is a path on which the stabilizers of vertices 
are ordered increasingly w.r.t.\ inclusion upon approaching $T'$. Obviously, the tree products of $T$ and $T'$ are isomorphic. 

\vv

{\bf (3.2) Example.} \ \ The mirror of a parabolic element can be contracted to 
a half-line.

\vv

{\bf (3.3) Definition.} \ \ Let $ T$ and $T'$ be subtrees of $\T_N^\ram$. We say that $T'$ is {\sl perpendicular} to $T$ (notation: ${T} \bot T'$)
if for every half-line $\ell = [\Lambda,\ast[ \ \subset T'$ and every edge $\sigma \in E(T)$ with
$\sigma \dashv \Lambda$, we have $\pi_N(\sigma) \neq \pi_N(\sigma')$, where $\sigma'$ is the unique edge $\sigma' \in \ell$ with $\sigma' \dashv \Lambda$. 

Since $\sigma$ and $\sigma'$ emanate from the same vertex, the last condition is equivalent to $\sigma \notin N_\Lambda \cdot \sigma'$. 

\vv

{\bf (3.4) Lemma. }\ \  {\sl Let $\Lambda \in V(\T_N^\ram)$. If there exists $\gamma \in N_\Lambda$ and $\sigma \dashv \Lambda$ such that $|\langle \gamma \rangle \cdot \sigma|>1$ {\rm (i.e.}, $\sigma$ is not fixed by all powers of $\gamma${\rm )}, then there exists a non-trivial $\gamma' \in N_\Lambda$ whose mirror is perpendicular to the $\gamma$-orbit of $\sigma$, {\rm i.e.}, such that $M(\gamma') \bot \langle \gamma \rangle \cdot \sigma.$}

\vv

\Proof Let $\nu=|\langle \gamma \rangle \cdot \sigma|$, and fix a numbering $\langle \gamma \rangle \cdot \sigma = \{ \sigma_i \}_{i=1}^\nu$. By the identification in (2.7), (letting $H=\rho(G), G=N_v$), it suffices to prove the following: if $p_1,...,p_{\nu+1}:=p_1$ is a sequence  of points on ${\bf P}^1$ which are cyclically permuted by an element $\gamma \in H$ (identified with $\sigma_1,...,\sigma_{\nu+1}=\sigma_1$), then  there exists a non-trivial $\gamma' \in H$ whose fixed points on ${\bf P}^1$ are not mapped to $\pi_H(p_1)(=\pi_H(p_i))$. 

Notice that $H \neq 1$ since $\gamma$ acts non-trivially, and that $\pi_H(p_1)$ is not totally ramified since $\nu >1$. It follows from (2.9) that if $H$ is cyclic, then every ramification point of $\pi_H$ must be totally ramified, and hence none of the $p_i$ are fixed by $H$. Hence we can let $\gamma'=\gamma$ in this case. 

If $H$ is not cyclic nor of type $B(t,n)$, then by (2.9), 
the map $\pi_H$ is branched above at least two points, and since (again, by (2.9)) there exist at least two points on ${\bf P}^1$ with disjoint decomposition groups, one can choose $\gamma'$ in the decomposition group of a ramified point that does not fix any point above $\pi_H(p_1)$. 

On the other hand, if $H=B(t,n)$, it suffices to take $\gamma' \in B(t,1) \subseteq B(t,n)$, since then, $\gamma'$ fixes a unique (totally ramified) point, which does not map to $\pi_H(p_1)$. \ee

\vv

{\bf (3.5) Lemma.} \ \ {\sl Let $\T^{\; \prime}$ be a subtree of $\T_N^\ram$ such that $\Ends(T_N^\ram) = \Ends(T)$, where $T=\pi_N(\T^{\; \prime})$. Then there is no half-line $\ell$ in $\T_N^\ram$ emanating from a vertex $\Lambda \in \T^{\; \prime}$ which is perpendicular to $\T^{\; \prime}$.}

\vv

\Proof  Suppose that such a half-line $\ell$ starting at
$\Lambda\in\T^{\; \prime}$ and pointing to $x \in {\bf P}^1(k)$ is perpendicular to $\T^{\; \prime}$. 
Then by our assumption on the ends, the image $\pi_N(\ell)$ is not a half-line, and hence there exists a vertex $\Lambda_0$  on $\ell$ which is not contained in $\T^{\; \prime}$, and $\gamma \in N_{\Lambda_0}$ which maps $[\Lambda_0,x[$ to a half-line containing $[\Lambda_0,\Lambda]$ (say we choose $\Lambda_0$ nearest to $\Lambda$). Set $\T_0:=\T^{\; \prime}\cup\ell$ and $v_0=\pi_N(\Lambda_0)$.
Let $T_0=\pi_N(\T_0)$, which is the union of
$T$ and the finite path $[v,v_0]$ (where
$v=\pi_N(\Lambda$)). Then:
\vv

(3.5.1) There exist an infinite sequence $\{\ell_n = [\Lambda_n,x_n[ \  \}_{n=0}^{\infty}$ of half-lines
in $\T^\ram_N$ and an infinite sequence $\{\varsigma_n\}_{n=0}^{\infty}$ 
of finite paths of length $>0$ in
$T^\ram_N$, such that for all $i,j$: 

(a) $\Lambda_0 \vdash \ell_0$ and
$T_0\cap\varsigma_0=\{v_0\}$,

(b) $\varsigma_i\subseteq\pi_N(\ell_i)$, 

(c) $\ell_i\cap\ell_j=\emptyset$ and
$\varsigma_i\cap\varsigma_j=\emptyset$ unless $|i-j|\leq 1$,

(d) $\ell_i\cap\ell_{i+1}=\Lambda_{i+1}$, 

(e) $\varsigma_i\cap\varsigma_{i+1}$ is a unique vertex, namely
$v_{i+1}:=\pi_N(\Lambda_{i+1});$ moreover,
$\varsigma_i=[v_i,v_{i+1}]$, 

(f) there exists an element $\gamma_i$ in
$N_{\Lambda_{i}}$ which maps $[\Lambda_i,x_{i-1}[$ to a half-line containing $[\Lambda_{i},\Lambda_{i-1}]$.

\unitlength1pt

\begin{center} \begin{picture}(90,90)
\put(10,35){\line(0,1){50}}
\multiput(10,33)(0,-2){3}{\circle*{1}}
\put(6,15){{\footnotesize $\T^{\; \prime}$}}
\put(80,80){{\footnotesize $\T_N^*$}}

\multiput(10,87)(0,2){3}{\circle*{1}}
\put(33,3){{\footnotesize $\ell_0$}}
\put(35,60){\vector(0,-1){50}}
\put(62,57){{\footnotesize $\ell$}}
\put(10,60){\vector(1,0){50}}
\put(87,32){{\footnotesize $\ell_1$}}
\put(35,35){\vector(1,0){50}}
\put(32,73){{\footnotesize $\gamma$}}

\put(35,65){\oval(10,10)[t]}
\put(30,62){\vector(0,-1){0}}
\put(40,60){\line(0,1){5}}
\put(16,33){{\footnotesize $\gamma_1$}}
\put(35,60){\circle*{2}}
\put(31,63){{\tiny $\Lambda_0$}}

\put(30,35){\oval(10,10)[l]}
\put(32,40){\vector(1,0){0}}
\put(30,30){\line(1,0){5}}
\put(56,48){{\footnotesize $\gamma_2$}}
\put(35,35){\circle*{2}}
\put(26,33){{\tiny $\Lambda_1$}}

\put(60,40){\oval(10,10)[t]}
\put(55,37){\vector(0,-1){0}}
\put(65,35){\line(0,1){5}}
\put(60,35){\circle*{2}}
\put(56,38){{\tiny $\Lambda_2$}}

\put(60,35){\line(0,-1){25}}
\multiput(60,8)(0,-2){3}{\circle*{1}}

\put(-30,-15){{\footnotesize Picture 2. The infinite sequence $\{ \ell_i \}$}}

\end{picture}
\end{center}

\vv

The proof of (3.5.1) is inductive. Since $\gamma$ does not fix all edges emanating from $\Lambda_0$, the conditions of (3.4) are satisfied. 
Hence there exists $\ell_0 \in \T^*_N$ starting at 
$\Lambda_0$ (say, pointing to $x_0 \in {\bf P}^1(k)$) which is perpendicular to $\T_0$.
Since the ends of $\pi_N(\T_0\cup\ell_0)$ coincide with 
those of $T$ by assumption, there has to be a vertex $\Lambda_1$ not in $\T_0$ and an element $\gamma_1 \in N_{\Lambda_1}$ which maps $[\Lambda_1,x_0[$ to a half-line containing $[\Lambda_1,\Lambda_0]$. Take such a vertex $\Lambda_1$
nearest to $\Lambda_0$, and set $v_1=\pi_N(\Lambda_1)$.
Let $\varsigma_0=[v_0,v_1]$.
Then (a) follows from the fact that $\ell_0$ is perpendicular to
$\T_0$, (b) is clear by definition, and we also have (f).

Assume that we have constructed all $\ell_j$ and $\varsigma_j$ for
$j<n$.
By lemma (3.4), there exists a half-line $\ell_n$ with $$ \ell_n \bot \T_n, \mbox{ where } \T_n:=\T_0\cup\,(\bigcup_{j<n}[\Lambda_j,\Lambda_{j+1}]).$$
Replacing $\T_0$ by
$\T_n$ in the argumentation above, we obtain $\Lambda_{n+1}$, $v_{n+1}$ and $\varsigma_n$.
Since a tree does not contain a cycle, all the properties (b)$\sim$(f)
are clear.
This finishes the proof of the (3.5.1).

Finally, 
$$\ell=[v,v_0]\cup\,(\bigcup_{n=0}^{\infty}\varsigma_n)$$ is a half-line in
$T^\ram_N$ such that 
$\ell\cap T=\{ v\}$. 
But this means that $T^\ram_N$ has an end which is not of
$T$, contradicting the assumption. \ee

\vv

(3.6) {\sl Proof of Proposition 1.} \ \ Fix a lift $\T^{\; \prime}$ of $T$ to $\T_N^*$. Let $\Lambda \in \T_N^\ram - \T^{\; \prime}$. Since $T_N^\ram$ is connected, there
exist a vertex $\Lambda' \in \T^{\; \prime}$ closest to $\Lambda$. The proof proceeds by induction on $d(\Lambda,\Lambda')$.

Assume $d(\Lambda,\Lambda')=1$. If $N_\Lambda=1$, then contraction is clear. Otherwise, take any non-trivial element $\gamma \in N_{\Lambda}$, and consider a half-line $\ell$ which starts at $\Lambda$ and converges to an end of
$M(\gamma)$.
We can then apply (3.5) to the tree $\T^{\; \prime \prime} := \T^{\; \prime} \cup [\Lambda,\Lambda']$ to find that $\ell$ is not perpendicular to $\T^{\; \prime \prime}$, and hence there exists
$\gamma'\in N_{\Lambda}$ such that $\gamma'\cdot\ell$ passes through 
the edge $\sigma=[\Lambda,\Lambda']$.
Since $\ell$ is pointwise fixed by $\gamma$, $\gamma'\cdot\ell$ is
pointwise fixed by $\gamma' \gamma {\gamma'}^{-1}$.
Hence $$ N_\Lambda = \bigcup_{\gamma' \in N_\Lambda} \gamma' N_\sigma {\gamma'}^{-1}, $$
{\sl i.e.,} $N_\Lambda$ is the union of the conjugates of its subgroup
$N_\sigma$. Hence (\cite{Huppert:79}, I.2, {\sl Aufgabe} 4) $N_\Lambda=N_\sigma$, so that
$\T^{\; \prime \prime}$ can be contracted to $\T^{\; \prime}$. It is clear how to proceed by induction. \ee

\vv

\vv

\section*{4. Mumford covers of ${\bf P}^1$ branched above at most three points} 

\vv

We will now describe the abstract group structure of $N$ if $\pi$ is only branched above $m=2$ or $m=3$ points. 

\vv

The possible structure of stabilizers of edges in $T_N$ turns out to
be further restricted by the following lemma:

\vv

{\bf (4.1) Lemma.} \ \ {\sl If $p>0$, then stabilizers of 
edges in $T_N$ are of the form $E_t \sd {\bf Z}_n$.}

\vv

\Proof As is shown in Proposition (4.5.1) of \cite{Gerritzen:80}, the stabilizer group of an edge contains a normal $p$-group with cyclic quotient. The statement follows from (2.9). \ee

\vv

We also have some information on the ``relative position'' of the $p$-prime
part of the stabilizer of an edge in the stabilizers of the corresponding vertices. We say
a subgroup of a group is {\sl maximally cyclic} if it is not contained in any bigger cyclic subgroup. 

\vv

{\bf (4.2) Lemma} (Herrlich \cite{Herrlich:80}, lemma 1). \ \ {\sl If it is non-trivial, then the prime-to-$p$ part of the stabilizer of an edge is a {\rm maximal} cyclic 
subgroup of the stabilizer of each of the corresponding vertices. \ee}

\vv

{\bf (4.3) Lemma.} \ \  {\sl Let $v \in V(T_N)$, $q=p^t$. If  $N_v = P(2,q)$, then 
for all $e \in E(T_N)$ with $e \dashv v$, either  $N_e = B(t, \{ \frac{1}{2} \}(q-1))$ or $N_e = {\bf Z}_{\{\frac{1}{2}\}(q+1)}$ {\rm (i.e.,} the $p$-part is also maximally cyclic in $N_v${\rm )}.}

\vv

\Proof  The action of $N_v$ on the edges emanating from a lift $\Lambda$ of $v$ to $\T_N$ induces a
representation $\rho : N_v \rightarrow PGL(2,\bar{k})$, which is faithful by (2.10). It suffices to invoke (2.7) and (2.9). \ee

\vv

{\bf (4.4) Lemma.} \ \ {\sl (i) The stabilizer of a vertex on the interior of a mirror of a parabolic element is a Borel group (interior meaning not an end point). 

(ii) If $\T_N^\ram$ contains a half-line $\ell$ which is pointwise fixed by a parabolic element of $N$, then $\pi_N(\ell) \subseteq T_N^\ram$ is not finite (whence, a half-line).}

\vv

\Proof (i) Let $\Lambda$ be such a vertex on $M(\gamma)$ for a parabolic
element $\gamma$ of $N$. By the description of $M(\gamma)$ given in (2.3), one sees that all
edges emanating from $\Lambda$ are fixed by $\gamma$; hence the group representation $\rho : N_\Lambda \rightarrow PGL(2,\bar{k})$ is not faithful, so the result follows from (2.10). 

(ii) Let $\ell$ be fixed by $\gamma$ and let $\Lambda$ be any interior vertex of $\ell \cap M(\gamma)$. It suffices to show that any two edges emanating from $\Lambda$ are mapped to different edges by $\pi_N$, {\sl viz.}, they are not identified by an element from $N_\Lambda$, which we know is a Borel group by (i). Let $e$ be
such an edge on $\ell$ pointing in the direction of the fixed point (say, $\infty$) of $\gamma$. Since all elements of $N_\Lambda$ fix the same point (they are upper-triagonal), in particular, $e$ is fixed by $N_\Lambda$, and hence cannot be identified with any other edge emanating from $\Lambda$. \ee

\vv

{\bf (4.5) Remark.} \ \ If the mirror of a parabolic element is contracted to
a half-line $[v_0,v_1,...[$ as in (3.2), its structure as a tree of groups is as follows: 

- the stabilizer $N_{v_0}$ can be larger than a Borel group (we see from (2.9) that is is then of type $P(2,q)$ since it should contain parabolic elements); 

- however, all interior
points are stabilized by a group $N_{v_i}=B(t_i,n)$ (this follows from part (i) of the above lemma), where $n$ is a fixed integer (by (4.2));

- the sequence $t_1,t_2,...$ is increasing since each $N_{v_i}$ stabilizes
the half-line $[v_i,...[$ (by the definition of mirror), and $N_{v_i} = N_{[v_i,v_{i+1}]} \subseteq N_{v_{i+1}}$.

\vv

{\bf (4.6) Proposition.} \ \ {\sl If $\pi : X \rightarrow S$ is branched above two
points, then $N$ is isomorphic to one of the following:

(A1)   $P(2,q)  \ast_{B(t,n_-)} B(t_1,n_-)$; 

(A2) $B(t_2,n_+) \ast_{{\bf Z}_{n_+}}   P(2,q)  \ast_{B(t,n_-)} B(t_1,n_-) $;

(A3) $B(t_3,n_-) \ast_{B(t,n_-)} P(2,q) \ast_{{\bf Z}_{n_+}}   P(2,q)\ast_{B(t,n_-)} B(t_1,n_-) $;

(A4) $B(t_4,n_+) \ast_{{\bf Z}_{n_+}}  P(2,q)$;

(A5) $B(t_5,n_-) \ast_{B(t,n_-)} P(2,q) \ast_{{\bf Z}_{n_+}}  P(2,q) $;

(B) $B(t'_1,n) *_{{\bf Z}_n} B(t'_2,n) $;

(C) $  E_{t'_3} * E_{t'_4}$, 

\noindent for $q=p^t, t|t_i (i=1,3,5), 2t|t_{2j} (j=1,2), n|p^{t'_1}-1,n|p^{t'_2}-1$ and $n_{\pm} = \{ \frac{1}{2} \} (q \pm 1)$.

If $p=3$, a list $(Ai')_{i=1}^5$ similar to $(Ai)_{i=1}^5$ can occur with $P(2,q),t,n_+, n_-$ replaced by $I,1,5,2$ respectively.

If $p=2$, there are the following extra possibilities:

(A1 $''$) $E_{t_2} \ast_{E_1} D_n \ast_{{\bf Z}_n} B(t_1,n)$;

(A2 $''$) $E_{t_2} \ast_{E_1} D_n \ast_{{\bf Z}_n} D_n \ast_{E_1} E_{t_1}$; 

(A3 $''$) $E_{t_2} \ast_{E_1} D_{n_+} \ast_{{\bf Z}_{n_+}} PGL(2,q) \ast_{B(t,n_-)} B(t_1,n_-), $

\noindent where $n$ is odd (and $n|p^{t_1}-1$ in the first case), $q$ is a power of $2$ and $n_\pm  =q \pm 1$.}

\vv

\Proof 
Let $T_0$ be a straight line in $T_N^\ram$ such that $\Ends(T_0)=\Ends(T_N^\ram)$; such a line exists, since we assume that there are only two branch points, so that $T_N^\ram$ has only two ends. By Proposition 1, we can contract 
$T_N^\ram$ to $T_0$ without affecting the structure of $N$. 

There exists at least one vertex in $T_0$ whose stabilizer is not classical, 
since there is wild ramification in $\pi : X \rightarrow S$. Let
$v$ be such a vertex. Then there are the following possibilities:

\vv

$(A)$ There is a vertex $v$ such that $N_v=P(2,q=p^t)$. By (4.4), the end emanating from a representative of $v$ in $\T_N^\ram$ which is
fixed by the Borel subgroup $B:=B(t,n_-)$ is mapped to a half-line of $T_0$ by $\pi_N$. Let $\{ v_i \}$ be a numbering of the vertices on this half-line in $T_0$ such that $v_i$ is at distance $i$ from $v=v_0$.
Since all $v_i$ are fixed by $B$, the stabilizer $N_{v_i}$ can only be a group of the form $B(t_i,n_i)$ by (2.9 iii). Indeed, if $N_{v_i} = P(2,p^{t_i})$, then $v_{i+1}$ would have to be stabilized by both ${\bf Z}_{n_+}$ and the Borel
group $N_{[v_{i-1}v_i]}=B(t_i,n_-)$ (using (4.3) and (4.4)), and hence by a group containing $B(2t_i,n_+)$, contradicting (4.3). 

Actually, $n_i = n_0=n_-$ by (4.2), and $\{t_i\}$ is increasing since $v_i$ is stabilized by the $N_{v_j}$ for $j \leq i$, and eventually terminates since $N_{v_j}$ converges to the (finite) stabilizer of the corresponding end (compare
with the argument in (4.5)).

From (2.7) and (2.9 iv), we see that the edge emanating from $v_0$ in the other direction is stabilized by ${\bf Z}_{n_+}$. For the stabilizer of the first neighbouring vertex $v'$ of $v$ on $T_0$ in the other direction for which $N_{v'} \neq N_v$, by (2.7) and (2.9) there are three possibilities: 

(A1) $N_{v'}={\bf Z}_{n_+}$; If all further vertices in the direction of
$v'$ have a cyclic stabilizer, they must all be equal by (4.2), and hence 
we are in case (A1). 

(A2-A5) On the other hand, should one of those
further vertices, say, $v"$, have a larger stabilizer, then we can first of all assume that $v"=v'$. Furthermore, if $p>3$, $N_{v'}$ can only be either (a) $N_{v'}=B(t',n_+)$ (with $q+1 | p^{t'}-1,t'>0$ by (4.2)) or (b) $N_{v'}=P(2,q)$ (with the same $q$ by (4.3)). 
Indeed, no other case from (2.9) can occur (dihedral or classical), because each time, lemma (3.4) applies to $N_{v'}$, leading to an extra
(third) end $M(\gamma')$ (in the notations of (3.4)) emanating from $v'$ in $T_0$, and thus contradicting (3.5). 

With only these two possibilities at hand, the subsequent neighbouring vertices of $v"$ in $T_0$ can be
dealt with as above, {\sl i.e.,} their stabilizers form an increasing stabilizing chain of Borel groups. This leads in the end to the prescribed form 
of $N$, which is (A2) in case (a) if $t_1>t$ and (A4) in case (a) if$t_1=t$; 
(A3) in case (b) if $t_1>t$ and (A5) in case (b) if $t_1=t$. 
\vv

$(A')$ Suppose we are not in case $(A)$, but $p=3$ and $N_v=I$ for some vertex $v \in T_0$. By (2.9), one end of $T_{N_v}^\ram$ is stabilized by the Borel group
$B(1,2)$ and one end by the cyclic group ${\bf Z}_5$. A similar reasoning as
in $(A)$ can be performed. 

\vv

$(B)$ Suppose we are not in case $(A)$ or $(A')$, but there is a vertex with $N_v=B(t,n),n>1$. One of the half-lines emanating from $v$ is stabilized by an increasing sequence of Borel groups as in case $(A)$. We can assume that the edge $e$ emanating from $v$ in the other direction is not fixed by parabolic elements (since all parabolic elements of $N_v$ have the same unique fixed point), whence $N_e = {\bf Z}_n$ by (4.2). Let $v'$ be the vertex closest to $v$ in 
this direction, such that $N_{v'} \neq {\bf Z}_n$ (which exists, since otherwise, $N$ would be finite). Then $N_{v'}$ can only be a Borel group (since it cannot be $P(2,q)$ by assumption, and it is not a classical group as in (A2-A5)). We find in the end that $N=B(t,n) *_{{\bf Z}_n} B(t',n)$. 

\vv

$(C)$ Suppose that all $N_v, v \in V(T_0)$ are $p$-groups, say, $N_v=E_{t_v}$. Because of lemma (4.4), there exist at least
two vertices $v,v'$ such that the mirrors of their respective stabilizers $N_v,N_{v'}$ map to different ends of $T_0$ under $\pi_N$. If we take such $v$ and $v'$ at minimal distance of each other, then all vertices between them are stabilized only by the trivial group, since the mirrors of elements in $N_v$ and
 in $N_{v'}$ have to map to disjoint half-lines in $T_0$. On edges outside of such a minimal segment, the $p$-groups are ordered increasingly and eventually stabilize. Hence, $N=E_t * E_{t'}$. 

\vv

The proofs if $p=2$ and there is a vertex $v$ on $T_0$ with $N_v=D_n$ are
entirely similar. \ee

\vv

{\bf (4.7)\ Proposition.} {\sl If $\pi : X\rightarrow S$ is branched above
three points, then $T_N^*$ can be contracted to a tree $T_0$ consisting of three
half-lines meeting at one common vertex $v_0$. Furthermore, either:

(D) There is an edge $e$ emanating from ${v_0}$ with $N_e=\{1\}$;

(E) The stabilizer of $v_0$ is classical non-dihedral, {\rm i.e.}, $N_{v_0} \in \{ T, O, I \}$; 

(F) $p\neq 2$ and $N$ is isomorphic to one of the following:
$$
\begin{picture}(270,50)(-160,-25)
\put(60,20){$Q(t_1)$}
\put(44,10){$\ast_{{\bf Z}_2}$}
\put(-130,0){$B(t_3,n_-)\ast_{B(t,n_-)}P(2,q)\ast_{{\bf Z}_{n_+}}D_{n_+}$}
\put(90,0){;}
\put(44,-12){$\ast_{{\bf Z}_2}$}
\put(60,-22){$Q(t_2)$}
\put(-180,0){(F1)}
\end{picture}
 $$
$$
\begin{picture}(270,50)(-160,-25)
\put(-43,20){$Q(t_1)$}
\put(-59,10){$\ast_{{\bf Z}_2}$}
\put(-130,0){$B(t,n)\ast_{{\bf Z}_n}D_n$}
\put(-13,0){,}
\put(-59,-10){$\ast_{{\bf Z}_2}$}
\put(-43,-20){$Q(t_2)$}
\put(-180,0){(F2)}
\end{picture}
 $$
where $q=p^t$, $t|t_3$, $Q(t_i)=B(t_i,2)$ $(i=1,2)$ and the two cyclic subgroups
${\bf Z}_2$ occuring in $(F)$ are the same in the dihedral group $D_n$. 

If $p=3$, a case $(F1')$ similar to $(F1)$ can occur, where $Q(t_i)$ is either $B(t_i,2)$ or
$PSL(2,3)\ast_{B(1,2)}B(t_i,2)$, and $P(2,q),t,n_+,n_-$ is
replaced by $I,1,5,2$, respectively.}

\vv

\Proof Let $T_0$ be the subtree of $T^{\ast}_N$ which is the union of three
half-lines $\ell_1$, $\ell_2$ and $\ell_3$ such that
$\ell_i\bigcap\ell_j=\{v_0\}$ for $i\neq j$ for some (uniquely determined) $v_0$. This $T_0$ exists, since we assume that only three points are branched in $\pi$, so that $T_N^\ram$ has three ends. 
By Proposition 1, we can contract $T^{\ast}_N$ to $T_0$.

Suppose $G=N_{v_0}$ is such that $\pi_G$ is branched over at most 2
points ({\sl cf.} (2.9)).
Then it follows from (2.7) that at least one of the edges emanating
from $v_0$ has trivial stabilizer; hence $(D)$ occurs.

Next, we suppose $\pi_G$ is branched over 3 points.
If we are not in case $(D)$, then $p\neq 2$ and the group $N_{v_0}$ is
a dihedral group $D_n$.
Let $v_1$ be the neighbouring vertex of $v_0$ such that the edge
connecting $v_0$ and $v_1$ is stabilized by ${\bf Z}_n$. For $N_{v_1}$, by
(2.7) and (2.9) there are three possibilities: (a) $N_{v_1}={\bf Z}_n$; (b)
$N_{v_1}=B(t,n)$ for $t>0$ or (c) $N_{v_1}=P(2,q)$.
Then an argument similar to that in the proof of (4.6) can be applied to
classify all the possible groups appearing in the amalgam product in
this direction.
The other directions and the case $(F1')$ can be treated similarly. \ee

\vv

\vv

\section*{5. The classical bound}

\vv

In this section, we will prove that $|\A(X)|$ is bounded ``classically'' by
$12(g-1)$ in case $(Ai')_{i=1}^5$, $(D)$ and $(E)$. The proof is very similar to Herrlich's computations in the $p$-adic case (\cite{Herrlich:80}). 

\vv

{\bf (5.1) Definition.} \ \ If $(T,G_\ast)$ is a finite tree of groups, we define its
$\mu$-invariant $\mu(T)$ to be
$$  (*) \ \ \mu(T) = \sum_{[vw] \in E(T)}  \frac{1}{|G_{vw}|} - \sum_{v \in V(T)} \frac{1}{|G_v|}.$$

\vv

The following group theoretical result provides the direct link between
$T_N$ as a tree of groups and the automorphism group of $X$.

\vv

{\bf (5.2) Theorem} (Karrass, Pietrowski, Solitar \cite{Karrass:73}). \ \ {\sl
With notations as above,
$ |N/\Gamma| \cdot \mu(T_N) = (g-1)$. \ee}

\vv

This theorem, together with (1.3), implies that in order to bound $|\A(X)|$ from above, one has to bound
$\mu(T_N)$ from below. 

\vv

{\bf (5.3) Lemma} (\cite{Gerritzen:80}, (4.7.1)). \ \ {\sl If in the above
situation, $B$ is a subtree of $T_N$,
then $\mu(T_N) \geq \mu(B)$. \ee}

\vv

{\bf (5.4) Proposition.} \ \ {\sl The bound $|\A(X)| \leq 12(g-1)$ holds in case $(D)$ and $(E)$ and $(Ai')_{i=1}^5$.}

\vv

\Proof Suppose that, in the notations of the proof of (4.6) and (4.7), $N_{v_0}$ is non-trivial such that $N_e=1$ for some
edge emanating from $v_0$ in $T_0$. Suppose that $e \in \ell_1$ for some half-line $\ell_1$ of $T_0$. Let $v_1$ be the vertex closest to $v_0$ on $\ell_1$ such that $N_{v_1} \neq 1$ (which exists, since eventually the stabilizer of $\ell_1$ is a (non-trivial) ramification group of $\pi$), and let
$B$ be the subtree of $T_0$ consisting of the path from $v_0$ to $v_1$. Then
$$\mu(B) = 1 - \frac{1}{|N_{v_0}|} - \frac{1}{|N_{v_1}|} \geq 1 - \frac{1}{2} - \frac{1}{3} = \frac{1}{6}, $$
unless $N_{v_0}$ and $N_{v_1}$ are both isomorphic to ${\bf Z}_2$ (whence, $p\neq 2$). Let $\Lambda_1$ be a lift of $v_1$ to $\T_N^\ram$. 
If $N_{v_1}={\bf Z}_2$ is generated by $\gamma$, then the mirror of $\gamma$
is an apartment $]x,y[$ in $\T_N$, and since $[v_0, v_1]$ is not stabilized by $\gamma$, by (3.4) there has to exist an element $\gamma' \neq \gamma \in N_{v_1}$ which maps $]x,\Lambda_1]$ to $[\Lambda_1,y[$; a contradiction. 

In case $(E)$, let $e$ be any edge emanating from $v_0$. We can assume
that $N_e \neq 1$, since otherwise, case $(D)$ occurs. Let $v_1$ be the vertex closest to $v_0$ on $T_0$ in the direction of $e$ such that $N_{v_1} \neq N_e$. Let $B$ be the segment connecting $v_0$ and $v_1$. We will now prove that
$\mu(B) \geq \frac{1}{12}$ for all such possible $B$. Recall that $N_e$ is maximally cyclic in both $N_{v_i} \ (i=0,1)$. 

\vv

(5.4.1) The following table
provides the necessary data in case $N_{v_1}$ is also a classical group:

\renewcommand{\arraystretch}{1.2}

{\footnotesize $$ \begin{tabular}{c|c||c|c}

 $N_{v_1},N_{e},N_{v_0}$ & $\mu_B$ & $N_{v_1},N_{e},N_{v_0}$ & $\mu_B$  \\

\hline

$D_{n_1}, {\bf Z}_{2}, T$ & $\geq \frac{1}{6}$ &

$D_{3}, {\bf Z}_{3}, T$ & $ =\frac{1}{12}$ \\

$D_{n_1}, {\bf Z}_{2}, O$ & $ \geq \frac{5}{24}$ &

$D_{3}, {\bf Z}_{3}, O$ & $ = \frac{1}{8}$ \\

$ D_{4}, {\bf Z}_{4}, O$ & $=\frac{1}{12}$ &

$D_{n_1}, {\bf Z}_2, I$ & $\geq \frac{7}{30}$ \\

$D_{3}, {\bf Z}_3, I$ & $= \frac{3}{20}$ &

$ D_5, {\bf Z}_5, I$ & $ = \frac{1}{12}$ \\

$T,{\bf Z}_{2,3},T$ & $=\frac{1}{3},\frac{1}{6}$ &

$T,{\bf Z}_{2,3},O$ & $=\frac{3}{8},\frac{5}{24}$ \\

$T,{\bf Z}_{2,3},I$ & $=\frac{2}{5},\frac{7}{30}$ &

$O,{\bf Z}_{2,3,4},O$ & $=\frac{5}{12},\frac{1}{4},\frac{1}{6}$ \\

$O,{\bf Z}_{2,3},I$ & $=\frac{53}{120},\frac{11}{40}$ &

$I,{\bf Z}_{2,3,5},I$ & $=\frac{7}{15},\frac{3}{10},\frac{1}{6}$ \\


\end{tabular} $$}

\renewcommand{\arraystretch}{1}

(5.4.2) Now suppose that $N_{v_1}$ contains a $p$-group. The classification theorem (2.9) shows that, except if $p=3$, we can assume that $N_e$ is cyclic of order $n$ prime to $p$; and it is maximally cyclic in both $N_{v_0},N_{v_1}$. This implies that $n=2,3$ if $N_{v_0}=T$, $n=2,3,4$ if $N_{v_0}=O$ and $n=2,3,5$ if $N_{v_0}=I$. If $N_{v_1}=PGL(2,q)$, then  $n=q + 1$ by (4.3). This leads to $q=3,4$ for $N_{v_0}=T$ and $q=2,3,4,5$ for $N_{v_0}=O$ and $N_{v_0}=I$. We can argue similarly using $2n=q + 1$ for $N_{v_1}=PSL(2,q)$. Consulting (2.9), the following possibilities remain:

\renewcommand{\arraystretch}{1.2}

{\footnotesize $$ \begin{tabular}{c|c}


$N_{v_1},N_{e},N_{v_0}$ & $\mu_B$  \\

\hline

$PSL(2,5),{\bf Z}_{3},T$ & $=\frac{7}{30}$ \\

$PSL(2,7),{\bf Z}_{4},O$ & $=\frac{17}{84}$ \\

$PSL(2,9),{\bf Z}_{5},I$ & $=\frac{13}{72}$ \\

$E_{t} \sd {\bf Z}_{2,3}, {\bf Z}_{2,3}, T$ &  $\geq \frac{19}{60},\frac{1}{8}$ \\

$E_{t} \sd {\bf Z}_{2,3,4}, {\bf Z}_{2,3,4}, O$ &   $\geq \frac{43}{120},\frac{41}{168},\frac{19}{120}$ \\

$E_{t} \sd {\bf Z}_{2,3,5}, {\bf Z}_{2,3,5}, I$ & $\geq \frac{23}{60},\frac{113}{420},\frac{109}{660}$ \\


\end{tabular} $$ }

\renewcommand{\arraystretch}{1}

(5.4.3) We are left to consider $p=3$ with $N_{e}=D_3 = E_1 \sd {\bf Z}_2$ and $N_{v_0}= I$. The following cases remain: (i) if $N_{v_1}=PGL(2,3)=O$, then $\mu(B)=\frac{13}{120}$; (ii) if $N_{v_1}=E_{t_1} \sd {\bf Z}_2$, then  $\mu(B)\geq \frac{17}{180}$; (iii) 
if $N_{v_1}=I$, then  $\mu(B)=\frac{2}{15}$. 

\vv

(5.4.4) In each of the cases $(Ai')_{i=1}^5$, one easily computes that
$\mu(T_0) \geq \frac{1}{12}$ from the definition (5.1). \ee

\vv

\vv

\section*{6. The extreme cases}

\vv 

We will now prove the bound in all remaining cases. Since the bound
is not linear in the genus, one cannot simply restrict to a suitably 
chosen segment of $T_0$ and apply the techniques of the previous section. 
Instead, we use the following lemma to further eliminate cases. 

\vv

{\bf (6.1) Notation.} \ \ Let $F$ be the real function
$$ F : {\bf R} \rightarrow {\bf R} \ : \ g \mapsto F(g):= 2 \sqrt{g} (\sqrt{g}+1)^2. $$
Let $p_\Gamma : N \rightarrow N / \Gamma$ be the projection map. 

\vv

Recall that for any finite group $G \subset N$, $p_\Gamma(G) \cong G$ since $\Gamma$ is free. 

\vv

{\bf (6.2) Lemma.} \ \ {\sl Let $X$ be a Mumford curve of genus $g$. Let $\lambda$ be the greatest common divisor
of $g-1$ and $|\A(X)|$, and write $g-1 = \lambda \cdot a$, $|\A(X)| = \lambda \cdot b$. If we can find a lower bound $\lambda_0$ for $\lambda$ such that 
$$ \lambda_0 b \leq F(\lambda_0 a + 1), $$
then the bound $|\A(X)| \leq F(g)$ holds.}

\vv

\Proof We have 
$$ |\A(X)| = (g-1) \frac{\lambda_0 b}{\lambda_0 a} \leq (g-1) \frac{F(\lambda_0 a +1)}{\lambda_0 a}. $$
Since the function $x \mapsto F(x)/(x-1)$ is increasing for $x \geq 2$ and
$\lambda_0 a +1 \leq \lambda a + 1 = g$, the result follows. \ee

\vv

{\bf (6.3) Proposition.} \ \ {\sl The bound $|\A(X)| \leq F(g)$ holds for all cases
$(Ai)_{i=1}^5$, except possibly $(A1)$ with $t_1=2t$ and $(A5)$ with $t_5=t$.}

\vv

\Proof We know that the ramification groups of $\pi \ : \ X \rightarrow S$ are
the stabilizers of the ends of $T_0$. 
The Riemann-Hurwitz formula applied to $\pi$ allows us to compute the ratio
of $g-1$ to $|\A(X)|$, and one 
finds the following values of $a$ and $b$ in each of the cases $(Ai)_{i=1}^5$:

\renewcommand{\arraystretch}{1.2}

{\footnotesize $$ \begin{tabular}{l|l|l|l}
& $a$ & $b$ & where \\
\hline
$(A1)$ & $q^n-q-1$ & $\{\frac{1}{2}\} q^n (q^2-1)$ & $nt={t_1}$ \\
$(A2)$ & $q^{n+1}-q^{n-m+1}-q^{n-m}-q+1$ & $\{\frac{1}{2}\} q^{n}(q^2-1)$ & $nt=2{t_2} \geq mt={t_1}$ \\
       & $q^{m+1}-q^{m-n+1}+ q^{m-n}-q-1$ & $\{\frac{1}{2}\} q^m(q^2-1)$ &
$nt=2{t_2} \leq mt={t_1}$ \\
$(A3)$ & $q^m-q^{m-n}-1$ & $\{\frac{1}{2}\} q^m(q-1)$ & $nt={t_1} \leq mt={t_3}$ \\
$(A4)$ & $q^{n+1}-q^{n}-q^{n-1}-q+1$ & $\{\frac{1}{2}\} q^{n}(q^2-1)$ & $nt={t_4}$ \\
$(A5)$ & $q^n-q^{n-1}-1$ & $\{\frac{1}{2}\} q^n(q-1)$ & $nt={t_5}$ \\
\end{tabular} $$}

\renewcommand{\arraystretch}{1}

\noindent One computes that the above lemma
 can be applied with $\lambda_0=1$ in all but the above mentioned cases. \ee

\vv

{\bf (6.4) Proposition.} \ \ {\sl The bound $|\A(X)| \leq F(g)$ holds in case $(A1)$ with $t_1=2t$, {\sl i.e.}, $n=2$ in the above table.}

\vv

\Proof Paradoxically, we will show that $|\A(X)|$ is bounded from {\sl below} non-trivially. 
Let $N_1 = P(2,q), N_2=B(2t,n_-), B=B(t,n_-)$ be the subgroups of $PGL(2,k)$ occurring in $T_0$, {\sl i.e.}, $N=N_1 \ast_B N_2$. We can assume that $N_1$ is the ``standard'' copy $P(2,{\bf F}_q)$ of $P(2,q)$ in $PGL(2,k)$ induced by the
inclusion ${\bf F}_q \subset k$.  Since the $p$-part of $N_1$ equals $B$, hence
is a subgroup of $N_2$, the $p$-parts of $N_1$ and $N_2$ have
to commute in $N$; so they can be simultaneously put into upper triagonal form, {\sl i.e.}, we can suppose that they are of the form
$$ \tau(\alpha_i) \in N_1 \cap N_2; \tau(\beta_j) \in N_2-N_1 $$ for some  $\{ \alpha_i \}_{i=1}^{t}, \ \{ \beta_j \}_{j=1}^{t} \in k$, where we set
$$ \tau(x) := \left( \begin{array}{cc} 1 & x \\ 0 & 1 \end{array} \right)$$
Now observe that the images of these matrices under the map $p_\Gamma\ :\ N \rightarrow N / \Gamma$ (where $\Gamma$ is the Schottky group of $X$), remain distinct. Indeed, if $p_\Gamma(\tau(x))=p_\Gamma(\tau(y))$, then $\tau(x-y) \in \Gamma$; but $\tau(x-y)$ is of finite order, so $x=y$ since $\Gamma$ 
is free. 

On the other hand, 
$N_1$ contains the ``non-split'' element $$\gamma = \left( \begin{array}{rr} 0 & 1 \\ -1 & 0 \end{array} \right),$$ so also the lower triagonal elements $\tau^t(x)=\gamma \tau(x) \gamma^{-1}$ are
contained in $N$ for all $x\in \{\alpha_i,\beta_j\}$. 
For the same reason as above, the elements $\tau^t(x)$ remain mutually distinct
modulo $\Gamma$. 

We will now prove that $p_\Gamma(\tau(x)) \neq p_\Gamma(\tau^t(y))$ for all $x,y \in \{ \alpha_i, \beta_j \}$. This is clear if $x,y$ are both in $ \{ \alpha_i \}$, since then
$\tau(x) \neq \tau^t(y)$ in  $N_1$ and $p_\Gamma(N_1)\cong N_1$. 

Now suppose that $p_\Gamma(\tau(\alpha_i)) = p_\Gamma(\tau^t(\beta_j))$ for some $i,j$. Since the commutator $[\alpha_i,\beta_j]=1$ in $N$, we have $$[p_\Gamma(\tau^t(\beta_i),p_\Gamma(\tau^t(\alpha_i))]=1.$$
We conclude that $$[p_\Gamma(\tau(\alpha_i)),p_\Gamma(\tau^t(\alpha_i))]=1.$$ But as $\tau(\alpha_i)$ and $\tau^t(\alpha_i)$ are both in $p_\Gamma(N_1) \cong N_1$, this gives a contradiction, since $\tau(\alpha_i)$ and $\tau^t(\alpha_i)$ are not commuting in 
$N_1=PGL(2,q)$.   

Finally, if $p_\Gamma(\tau(\beta_i))=p_\Gamma(\tau^t(\beta_j))$, then 
since the first element commutes with $p_\Gamma(\tau(\alpha_1))$ and the second with 
$p_\Gamma(\tau^t(\alpha_1))$, we would find that $p_\Gamma(\tau(\alpha_1))$ and
$p_\Gamma(\tau^t(\alpha_1))$ commute, leading to the same contradiction. 

We find in the end that the order of $A=N/ \Gamma$ is at least divisible by $q^4$ ({\sl viz.}, the number of distinct $p$-order elements constructed above), which means that we can choose $\lambda_0=q^2$ in 
lemma (6.2) and this leads to the desired bound, except if
$q=2$, but then, $PGL(2,2)=D_3$ and another case applies. \ee

\vv

{\bf (6.5) Proposition.} \ \ {\sl The bound $|\A(X)| \leq F(g)$ holds in case $(A5)$ with $t_5=t$.}

\vv

\Proof We see that the group $p_\Gamma({\bf Z}_{n_+})={\bf Z}_{n_+}$ which stabilizes the central edge of $T_0$ acts
without fixed point on $X$ (by the correspondence between fixed points and ends). By Hurwitz's formula applied to this unramified action, $n_+$ divides $g-1$. Hence we can set $\lambda_0=n_+$ to see that the criterion of lemma (6.2) is satisfied. \ee

\vv

{\bf (6.6) Lemma.} \ \ (\cite{Serre:80}, I.1., Prop.\ 4) \ \ {\sl If $G_1$ and
$G_2$ are finite groups, then we have a natural exact sequence
$$ 1 \rightarrow [G_1,G_2] \rightarrow  G_1 * G_2 \rightarrow G_1 \times G_2 \rightarrow 1, $$ where the group $[G_1,G_2]$ generated by the
commutators $[g_1,g_2]:=g_1 g_2 g_1^{-1} g_2^{-1}$ for $g_1 \in G_1-\{1\}$ and $g_2 \in G_2-\{1\}$ is 
a maximal free subgroup of $G_1 * G_2$ of rank $(|G_1|-1)(|G_2|-1).$ \ee}

\vv

{\bf (6.7) Lemma.} \ \ {\sl Let $X$ be an ordinary curve, admitting a Galois
cover $X \rightarrow {\bf P}^1$ branched above two points. If its Galois group is of the form $G=B(t,n)$, then $n=2$.}

\vv

\Proof  Let $\bar{E}$ be a subgroup of 
index $p$ in a $p$-Sylow group of $G$ and consider the curve $Y:=\bar{E} \backslash X$, which is 
a Galois cover of $G \backslash X = {\bf P}^1$, totally ramified above two points with 
Galois group ${\bf Z}_p \sd {\bf Z}_n$. It follows that the intermediate 
curve $Y':= {\bf Z}_p \backslash Y$ is a cyclic $p$-prime Galois cover of ${\bf P}^1$ ramified above two points, and hence equal to ${\bf P}^1$. Let $x$ be a coordinate on $Y'$; then
the equation of $Y$ is of the form
$$ y^p-y = B, \ B := c \frac{x-\alpha}{x} \mbox{ \ \  for some } \alpha \neq 0 $$
and $c \in k$. Indeed, the ramification points are the poles of $B$, and since $X$ is ordinary, they are simple (Subrao \cite{Subrao:75}). This curve has to admit an automorphism of order $n$. It is easy to see that
such an automorphism has to be of the form $ B \mapsto \zeta B + B_0^p-B_0, y \mapsto \zeta y + B_0$ for some $n$-th root of unity $\zeta$ and some $B_0(x) \in k(x)$ ({\sl cf.} Hasse \cite{Hasse:35}, p.\ 38).  On the other hand, the automorphism should be induced from an automorphism of $k(x)$ of order $n$, {\sl i.e.}, of the form $x \mapsto \zeta \cdot x$. This leads to the identity
$$ B(\zeta x ) - \zeta B(x) = B_0(x)^p-B_0(x), $$
which is only satisfied if $\zeta=\pm 1$ with $B_0$ a (constant) root of $c(1-\zeta)=B_0^p-B_0$. Hence $n=2$ as claimed. \ee

\vv

{\bf (6.8) Proposition.} \ \ {\sl The bound $|\A(X)| \leq F(g)$ holds in case $(B)$ and $(C)$.}

\vv

\Proof Recall that in this case, $N=N_1 \ast_B N_2$ with $N_1=B(t,n)$, $N_2=B(t',n)$ and $B={\bf Z}_n$ for $n|(p^t-1;p^{t'}-1)$. The criterion in 
lemma (6.2) fails to hold with $\lambda_0=1$ only if $t=t'$, which we will 
assume from now on. Observe that then, $N = (E_1 \ast E_2) \sd {\bf Z}_n$, where $E_1$ and $E_2$
are the respective $p$-parts of $N_1$ and $N_2$.  

We will distinguish two cases. First of all, assume that $E_1=E_2 \mbox{ mod } \Gamma$. It follows from this that $[E_1 ,E_2 ] \subseteq \Gamma$. 
In particular, we find that $\A(X)=N/\Gamma$ is a quotient of $(E_1 \ast E_2) \sd {\bf Z}_n / [E_1 ,E_2 ] = (E_1 \times E_2) \sd {\bf Z}_n$ by (6.6). Since $E_1=E_2 \mbox{ mod }
\Gamma$, $\A(X)$ is even a quotient of $E_1 \sd {\bf Z}_n$. On the other hand, 
since $E_1 \sd {\bf Z}_n \cong N_1 \subseteq N/\Gamma$, we find in the end that
$$ \A(X) = E_1 \sd {\bf Z}_n. $$
The previous lemma implies that $n=2$, and knowing this, the bound $\mu(T_0) \geq \frac{1}{12}$ is easy to check. 

Let $E_1'$ denote the part of $E_1$ which is identified with a part (called $E'_2$) of $E_2$ modulo $\Gamma$. Denote by $E_i''$ the respective complements
of $E_i'$ in $E_i$, and let $p^\tau$ be the order of $E_i''$. By the result of the previous
paragraph, we can assume that $\tau>0$. We claim that $E_i''$ are
${\bf Z}_n$-modules (for the action of conjugation). Indeed, it suffices to show
this for $E_i'$. Let $\varepsilon$ be an element of $E_1'$, and choose an
element $\bar{\varepsilon} \in E_2'$ with $p_\Gamma(\varepsilon) = p_\Gamma(\bar{\varepsilon})$. We have to show that for every $\sigma \in {\bf Z}_n$, $\sigma \varepsilon \sigma^{-1} \in E_1'$, {\sl viz.}, $p_\Gamma(\sigma \varepsilon \sigma^{-1}) = p_\Gamma(x)$ for some $x \in E_2$, and it suffices
to let $x=\sigma \bar{\varepsilon} \sigma^{-1}$. 

Since ${\bf P}^1$ does not admit unramified extensions, the automorphism
group of $X$ is generated by its decomposition groups, and in particular, 
it contains the group generated by $E_1 \sd {\bf Z}_n$ and $E_2 \sd {\bf Z}_n$. By assumption, $E_1$ and $E_2$ intersect modulo $\Gamma$ in a group of 
order $p^{t-\tau}$, so that 
$$ |\A(X)| \geq \frac{|E_1||E_2||{\bf Z}_n|}{|E_1 \cap E_2|} = p^{t+\tau}n,$$
hence we can let $\lambda_0=p^\tau$ in (6.2). If we then use that $n \leq p^\tau-1$ (since $E_1''$ is a ${\bf Z}_n$-module, $n$ divides $p^\tau-1$), we find that the inequality in (6.2) is fulfilled, so that the desired bound holds (note that if $p^t=2$, we have $g=1$). \ee

\vv

{\bf (6.9) Proposition.} \ \  {\sl The bound $|\A(X)| \leq F(g)$ holds in case $(F1), (F1')$ and $(F2)$.}

\vv

\Proof In all such cases, one can compute a corresponding table
of $(a,b)$ depending on the values of $t,t_1,t_2$. We will skip the 
detailed computation; let it suffice to say that lemma (6.2) cannot be 
applied with $\lambda_0=1$ only if $t_1=t_2=0$, in which case the order 
of the ramification groups is $(2,2,nq)$ (with $q=p^t$) and $a=q-2, b=2nq$
(respectively, $a=\frac{1}{2}(q-2), b=nq$ if $p=2$). 

For the case $(F1)$, observe that, as in the proof of (6.5), the group ${\bf Z}_{n_+}$ acts freely on $X$, leading to $\lambda_0 \geq n_+$, and with this 
value, the bound follows from (6.2). 

The only case that remains to be settled is $(F2)$, {\sl viz.}, $N=N_1 \ast_B N_2$, with $N_1=B(t,n), B={\bf Z}_n$ and $N_2= D_n$. Let $E$ be a $p$-Sylow of $N_1$, let $\gamma$ be an involution in $D_n-{\bf Z}_n$ and set $$E'=\{ \varepsilon \in E : p_\Gamma(\gamma \varepsilon \gamma) \in p_\Gamma(E) \}. $$ 
By lemma (6.7), if $E=E'$, then $n=2$ and the bound follows easily. 

Let $E''$ be the complementary ${\bf Z}_n$-module of $E'$ in $E$, and let 
$|E''|=p^\tau$; then we can assume that $\tau>0$. We find as in (6.8) that
$\A(X) \geq 2 np^{t+\tau}$, and the desired bound follows from (6.2) taking
$\lambda_0= 2 p^\tau$. \ee

\vv

{\bf (6.10) Remark.} \ \ Stichtenoth has shown in \cite{Stichtenoth:73} that $\A(X) \leq 84(g-1)$ if $\pi$ is ramified above three points and at least two ramification indices are $>2$. 

\vv

{\bf (6.11) Proposition.} \ \ {\sl The bound $|\A(X)| \leq F(g)$ holds if $p=2$ in case $(Ai'')_{i=1}^3$.}

\vv

\Proof Unconditionally in case $(A2'')$ and in case $(A1'')$ and $(A3'')$ if $t_2 \neq 1$, it is easy to compute that $\mu(T_0) \geq \frac{1}{12}$. 

In case $(A1'')$ with $t_2=1$, one finds that $a=2^{t_1-1}-1,b=2^{t_1}n$, and $\pi$ has ramification type $(2,2^{t_1}n)$. As in (6.9), let $E$ be the $p$-part of $B(t_1,n)$ and $E'$ the image of the conjugate of $E$ by an involution of $D_n$. If $E=E'$, then again, $n=2$ and the bound holds. If we let $\tau$ be the order of a complement of $E'$ in $E$, then we can set $a=2^{t_1}-1,b=2^{t_1}n,\lambda_0=2^{\tau+1}$ and using $n \leq 2^\tau-1$, prove that (6.2) is satisfied (note that $t \geq \tau$).  

In case $(A3'')$ with $t_2=1$, we define $n$ by $tn=t_1$; then $a=\frac{1}{2}q^n-1, b=q^n(q-1)$, and since ${\bf Z}_{n_+}$ acts freely on $X$, we can let $\lambda_0=q+1$ as in (6.5). \ee 

\vv

{\bf (6.12) Remark.} \ \ A careful inspection of the estimates shows that, if $p \neq 2$, only
case $(F2)$ with $t_1=t_2=0, n=p^t-1$ exactly attains the bound of the main theorem; similarly, if $p=2$, then only case $(A1'')$ with $t_1=t,t_2=1,n=p^t-1$ does so (if a free subgroup of $N$ of rank $g=(p^t-1)^2$ with normalizer $N$ exists, {\sl cf. infra}). Observe that in both cases, the group $N$
has the same form (amalgamation of a dihedral with a Borel group over a cyclic
group), but there are 3 branch points if $p \neq 2$, whereas if $p=2$, there 
are only 2 branch
points (compare with the two different cases in (2.9ii)). 

\vv

\vv

\section*{7. The icosahedral group -- proof of proposition 2}

\vv

(7.1) We have already remarked that $12(g-1)>F(g)$ only if $g \in \{5,6,7,8\}$, 
and hence there is only a finite number of groups $A$ (namely, 134) for which 
$12(g-1) < |A| \leq F(g)$. These, one can easily write down, {\sl e.g.}, using
{\tt GAP} (\cite{GAP}). The only non-solvable such group is $I$, whose order
is 60. It follows that $g=6$. Note that $I$ is a subgroup of $PGL(2,k)$ only if $p \neq 2,5$, and then, a geometric construction of a Mumford curve $X$ of genus 6 with 60 automorphisms was given by Herrlich in \cite{Herrlich:78}, {\sl pp.} 50-51. The normalizer of its Schottky group is of the form
$$ N = I \ast_{{\bf Z}_5} D_5, $$
and $X$ is a cover of ${\bf P}^1$ ramified above 4 points with ramification 
indices $(2,2,2,3)$ if $p\neq3$ and $(2,2,6)$ if $p=3$. 

\vv

(7.2) We will now compute the dimension of the stratum of the moduli space
$M_6$ (which we consider as an algebraic space) of curves of genus 6 containing this particular example. Observe that
the dimension of the infinitesimal deformation space of $X$ consisting of curves whose automorphism group
contains $I$ is given by 
$ h^1(X,T_X)^I$, where $T_X$ is the tangent sheaf of $X$ and the superscript
$I$ means taking group invariants. By Serre duality, 
$ h^1(X,T_X)^I = h^0(X,\Omega^{\otimes 2}_X)^I, $
where $\Omega_X$ is the sheaf of regular differentials on $X$. Finally, this space of invariant 2-differentials is computed for $(p,|I|)=1$ to be of dimension
$3 \bar{g} - 3 + n$, where $\bar{g}$ is 
the genus of $I \backslash X$ and $n$ is the number of branch points of $X \rightarrow I \backslash X$ (\cite{Farkas:80}, {\sl pp.} 254-255). Hence in our case, the closed stratum $\bar{S}_I$ of curves in $M_6$ whose automorphism group contains $I$ is one-dimensional. We do not know whether $\bar{S}_I$ is connected, and whether or not it equals the open stratum $S_I$. 

\vv

(7.3) Let us now switch to ${\cal M}_g$, the moduli space of Mumford curves of genus g over $k$ (which can conveniently be described by non-archimedean Teichm{\"u}ller theory, {\sl cf.} Herrlich \cite{Herrlich:84}).
Recall that $K$ is the algebraic closure of $k$. We use the following (folklore) fact about this moduli 
space:

\vv

{\bf (7.3.1) Claim.} \ \ {\sl ${\cal M}_{g,K}$ is an open analytic subspace of the analytification $M^{\an}_{g,K}$ of the moduli space $M_{g,K}$.} 

\vv

A sketch of the proof goes as follows: For each of the moduli spaces $M$ that we consider, let $\bar{M}$ denote its Deligne-Mumford compactification. First one resolves the technical 
problem that $\bar{M}_{g,K}$ is not a {\sl fine} moduli space by adding some
level-$n$ structure for $n$ large enough so that $\bar{M}_{g,K}[n]$ is a fine moduli space, and note that $\bar{M}^{\an}_{g,K}[n] \rightarrow \bar{M}^{\an}_{g,K}$ is 
an open map with respect to the analytic topology. This follows because (by a short local calculation using the fact that the corresponding algebraic map $\bar{M}_{g,K}[n] \rightarrow \bar{M}_{g,K}$is algebraically finite flat), it is analytically flat quasi-finite in the sense
of Berkovich (\cite{Berkovich:93}, 3.2), and applying (3.2.7) of {\sl loc.\ cit.} 

Let $\bar{M}^{\an}_{g,K}[n] \rightarrow \bar{M}^{\an}_{g,\bar{K}}[n]$ be the reduction map that
associates to any curve its stable reduction over the residue field $\bar{K}$ of $K$, and let $Z$ be the locus of multiplicative reduction on $\bar{M}^{\an}_{g,\bar{K}}[n]$. Then, by deformation theory, $Z$
is a Zariski closed set, say, given locally by 
equations $f_1=...=f_r=0$. The locus of Mumford curves 
in $\bar{M}^{\an}_{g,K}[n]$ is given by the ``tube'' $]Z[$ of $Z$, where $]Z[$ is defined locally by
$$ ]Z[ :=\{ x \in \bar{M}^{\an}_{g,K}[n] \ : \ |F_i(x)|<1 \mbox{ for all } i=1,...,r \}, $$
 where $F_i$ are lifts of $f_i$ to $\bar{M}^{\an}_{g,K}[n]$. Since $]Z[ = \lim\limits_{\stackrel{\longrightarrow}{e}} \  ]Z[_e$ where the limit is over all ``values'' $e \in |K^*|$, and 
$$ ]Z[_e = \{ x \in \bar{M}^{\an}_{g,K}[n] \ : \ |F_i(x)| \leq e \mbox{ for all } i=1,...,r \} $$
 are affinoid open sets, we find that $]Z[$ is open with respect to the strong rigid topology. The claim follows. \ee

\vv

It follows immediately from the computation in (7.2) that the stratum $\bar{S}_I$ intersects ${\cal M}_6$ in a one-dimensional space. Since we know that $|\A(X)| \leq 60$ on ${\cal M}_6$, we have
$S_I \cap {\cal M}_6 = \bar{S}_I \cap {\cal M}_6$, hence all curves in this one-dimensional rigid analytic stratum have automorphism group exactly equal to $I$.

However, there is a more direct way to the computation of the dimension of closed strata in the space of Mumford curves that avoids the use of invariant theory (hence, of the assumption that the characteristic is coprime to the order of the automorphism group):

\vv

{\bf (7.4) Theorem} (Herrlich \cite{Herrlich:84}). \ \ {\sl With the above notations, let $A=N/\Gamma$. If $g \geq 4$, then the closed stratum $\bar{S}_A \subset M_g$ of curves whose automorphism group contains $A$ intersects the Mumford locus ${\cal M}_g$ in a space of dimension
$$ \dim(\bar{S}_A \cap {\cal M}_g) = 3(f+d_v-d_e-1)+2(c_v-c_e), $$
where $c_v$ (resp.\ $d_v$) is the number of non-trivial cyclic (resp.\ non-cyclic) vertex groups of $T_N$, $c_e,d_e$ are similarly defined for edges of $T_N$, and $f$ is the number of free generators of $N$.} 

\vv

\noindent Using (7.4), one computes that in our case,  $\dim(\bar{S}_I \cap {\cal M}_6)=1$ as expected, and the same result as in (7.3) follows. Note, however, that (7.4) does not say anything about the dimension of the corresponding stratum in $M_g$.

\vv

{\bf (7.5) Remark.} \ \ Dolgachev has pointed out to the authors that the one-parameter family of curves $S_I$ can be made explicit by a classical geometric construction as a pencil of curves on the Del Pezzo quintic, {\sl cf.} Edge \cite{Edge:81}. 

\vv

{\bf (7.6) Remark.} \ \ The original version of theorem (7.4) also included
the statement that the open stratum $S_A \cap {\cal M}_g$ is rigid-analytically {\sl connected}. However, a careful study revealed that this connectedness statement, which is theorem 2 in \cite{Herrlich:84}, is violated by the pencil mentioned in (7.5). Herrlich agrees with the authors of the current paper that this is due to the fact that proposition 1 in \cite{Herrlich:84} is not true as it is stated. However, this does not further affect the results of this paper. A more detailed analysis of the connectedness of the stratification in the Mumford locus can be found in \cite{Kato:00}.

\vv

\vv

\section*{8. Discreteness -- attaining the bound}

\vv

(8.1) It has not yet been proven that the groups on the  lists in (4.6) and (4.7) actually occur as discrete subgroups of $PGL(2,k)$. It happens exactly when there is an action of $N$ on $\T$ whose quotient equals $T_N$. For this, it suffices that the action of $N$ is well-defined, and has a finite stabilizer at at least one point of $\T$. Checking this is not always so
obvious, but there is a direct criterion to see whether a {\sl free} product exists as a discrete subgroup of $PGL(2,k)$. 

\vv

{\bf (8.2) Definition.} The {\sl isometric circle} of a non-trivial element $\gamma$ of 
finite order in $PGL(2,k)$ is defined by
$\bar{I}_\gamma = \{ z \in {\bf P}^1(k) : |cz+d | \leq 1 \}, $ where $(c,d)$ is the
second row of $\gamma$. 

\vv

{\bf (8.3) Lemma} (\cite{Herrlich:80}). \ \ {\sl If $G,H$ are finite subgroups of $PGL(2,k)$, then $G \ast H$ exists as a discrete subgroup of $PGL(2,k)$ if and only if $\bar{I}_\gamma \cap \bar{I}_\delta = \emptyset$ for all non-trivial $\gamma \in G$ and $\delta \in H$.}

\vv

(8.4) We will use the above criterion to prove that the segment that attains
the bound of the main theorem ($(F2)$ with $t_1=t_2=0, n=p^t-1$ if $p \neq 2$  and $(A1'')$ with $t_1=t,t_2=0,n=p^t-1$ if $p=2$) does correspond to a discontinuous group, and hence a corresponding Mumford curve exists. Let us consider the 
segment $\B$ with $N_1 = B(t,n)$ and $N_2=D_n = {\bf Z}_n \sd \langle \gamma \rangle$ for some element $\gamma$ of order 2. We will study this case 
as a kind of ``cover'' of case $(B)$. 
Let $\B'$ be the segment defined by $N_1=B(t,n) = E_t \sd {\bf Z}_n$, $N_2=E'_{t} \sd {\bf Z}_n$ and $H={\bf Z}_n$, where $E_{t},E'_{t}$ are two disjoint $p$-groups generated by $\{ \varepsilon_i \}_{i=1}^{t_1}$ and $\{ \varepsilon'_i \}_{i=1}^{t_1}$ respectively, and on which, by definition, ${\bf Z}_n$ acts componentwise by conjugation. Looking at explicit presentations, one sees that 
$$ N_{\B} = N_{\B'} * \langle \gamma \rangle / \langle \gamma^2 = 1, 
\gamma \sigma \gamma = \sigma^{-1}, \gamma \varepsilon_i \gamma = {\varepsilon_i'}, \ i=1...t \rangle, $$ 
with obvious notations.

\vv

{\bf Remark.} \ \ Geometrically, the link between the segments $\B$ and $\B'$ is
as follows: one subdivides the segment $\B$ like $$ (E_{t} \sd {\bf Z}_n) *_{{\bf Z}_n} {\bf Z}_n *_{{\bf Z}_n} (E'_{t} \sd {\bf Z}_n),$$ and then lets an element
$\gamma$ of order two act on this as a vertical ``mirror'' along the central vertex; via $\gamma$, one identifies $E_{t}$ and $E'_{t}$ by conjugation and the elements on the stabilizers of edges by inversion. 

\vv

Rearranging words shows that $N_{\B} = \langle \gamma \rangle \cdot N_{\B'}$, so that $N_{\B'}$ is a subgroup of (finite) index 2 in $N_{\B}$ (of course, we actually mean the image of $N_{\B'}$ in $N_\B$). To show
that $N_{\B}$ is discontinuous, it therefore suffices to show this for $N_{\B'}$.  One sees that 
$ N_{\B'} = (E_{t} * E'_{t}) \sd {\bf Z}_n$.
We can let $E'_{t} = Q E_{t} Q^{-1}$ for some $Q \in PGL(2,k)$.  It suffices to show that the free product $E_{t} * E'_{t}$ is contained 
in $PGL(2,k)$. Now the isometric circles of elements in $E'_{t}$ are translates by $Q$ of the corresponding
circles in $E_{t}$, which we can assume to be unit circles with centers at finitely many
elements of $\bar{{\bf F}}_p$ (if we let $E_{t}$ be generated by suitable lower triagonal matrices). Call their union $I$. It then suffices to choose $Q$ such that $Q \cdot I \cap I = \emptyset,$ which is clearly possible since $k$ contains 
elements of non-unit valuation, so that we can let $Q$ be translation over an element of high enough valuation.

The corresponding Schottky group $\Gamma_t = [E_t, \gamma E_t \gamma]$ is free of rank $g=(p^t-1)$, and normal in $N$. By our main theorem, $N$ is then exactly
equal to the normalizer of $\Gamma$. 

\vv

\section*{9. Artin-Schreier-Mumford curves -- proof of proposition 3} 

\vv

(9.1) Let $q=p^t$. To avoid problems with the singularities of the particular plane model given, we think of the curve $X_{t,c} : (y^q-y)(x^q-x)=c$ as embedded in ${\bf P}^1 \times {\bf P}^1$. Since $|c|<1$, its analytic reduction is given by the ``chess board'' of lines $(x^q-x)(y^q-y)=0$, and hence it is a Mumford curve whose genus (the number of squares of the board) is given by $g_t=(q-1)^2$.  It seems appropriate to call these curves {\sl Artin-Schreier-Mumford} curves. 
\unitlength0.7pt
\begin{center}
\begin{picture}(110,120)
\put(47,110){{\footnotesize $q$}}
\put(10,95){$\overbrace{\mbox{\hspace{2cm}}}$}
\put(-32,42){{\footnotesize $ q \left\{ \begin{array}{c} \\ \; \\ \; \\ \; \\ \; \\  \end{array} \right. $}}
\multiput(10,0)(20,0){2}{\line(0,1){90}}
\multiput(92,0)(20,0){1}{\line(0,1){90}}
\multiput(0,7)(0,20){1}{\line(1,0){100}} 
\multiput(0,63)(0,20){2}{\line(1,0){100}}
\multiput(42,48)(6,-4.7){7}{\circle*{1,4}}
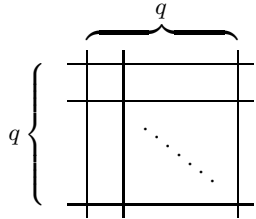
\put(-80,-25){{\footnotesize Picture 3. the analytic reduction of $X_{t,c}$}}
\end{picture}
\end{center}

\vv

\vv

(9.2) One can see from the equation that the curve admits the following automorphisms:
$(x \mapsto \alpha x, y \mapsto \alpha^{-1} y)$ for $\alpha \in {\bf F}_q^*$, 
$(x \mapsto x+\beta, y \mapsto y + \gamma)$ for $\beta, \gamma \in {\bf F}_q$ and
$(x \mapsto y, y \mapsto x)$. This accounts for $2q^2(q-1)=F(g_t)$ automorphisms, 
and hence by the main theorem, there are no more, and $N_t$ is as expected (note that this gives a less intricate proof of the fact that $N_t$ is discrete). It is also clear that the 
automorphism group is given by
$ \A(X_t) = {\bf Z}_p^{2t} \sd D_{p^t-1}, $
and hence the Schottky group (which is the kernel of $N_t \rightarrow \A(X_t)$) is generated by the commutators $[{\bf Z}_p^t, \gamma {\bf Z}_p^t \gamma]$ for
a fixed involution $\gamma \in D_{p^t-1}-{\bf Z}_{p^t-1}$. 

\vv

(9.3) The dimension of the stratum of curves $X$ in $M_{g_t}$ whose automorphism group contains $\A(X_{t,c})$ can not be determined in a straightforward way because 
of wild ramification in $\pi$ (\cite{Farkas:80} does not apply immediately), but Herrlich's formula (7.4) does imply that the intersection of the closed stratum 
$\bar{S}_{\mbox{{\tiny Aut}}(X_t)}$ with the Mumford locus is a one-dimensional subspace. By the main theorem, this is also the intersection with the open stratum. Note that the family $X_{t,c} \rightarrow {\bf A}^1_k$ (given by projection onto $c$) is not constant as $c$ varies, since it is a stable curve having a singular fibre over the origin (so it has moving moduli). As $c$ varies through $\{ |c|<1 \}$, we move through the open stratum in ${\cal M}_{g_t}$ (and  $\gamma$ changes). 

\vv

\vv

\section*{10. Drinfeld modular curves -- proof of proposition 4}

\vv

{\bf (10.1) Lemma} (Gekeler \cite{Gekeler:86}, VII.5). \ \ {\sl The cover $X({\frak n}) \rightarrow X(1) = {\bf P}^1$ of Drinfeld modular curves is ramified above two points with
respective indices $q+1$ and $q^d(q-1)$, where $d=\deg({\frak n})$. The genus of 
$X({\frak n})$ satisfies
$$ g(X({\frak n}))-1 = |G({\frak n})| \frac{q^d-q-1}{q^d(q^2-1)}, \mbox{ where } G({\frak n}) = \Gamma(1)/\Gamma({\frak n}){\cal Z} $$ is of order
$$ |G({\frak n})| = q^{3d} \prod_{p|n} (1 - \frac{1}{q^{2\deg(p)}}). $$}

\vv

{\bf (10.2) Notation.}  If $G$ is a group of automorphisms
acting on a curve $X$, we denote by $G_{i,P}$ for a point $P \in G \backslash X$ and $i \in {\bf Z}_{\geq 0}$ the $i$-th ramification group
of $X \rightarrow G \backslash X$ at $P$. As usual, let $\pi_G$ denote the corresponding 
covering morphism.

\vv

The following lemma characterizes (in a very special situation) when a group of automorphisms of a curve is
the full group of automorphisms in terms of ramification data. 

\vv

{\bf (10.3) Lemma.} \ \ {\sl Assume $p \neq 2$. Let $A$ be a group of automorphisms acting on an ordinary curve $X$ with quotient $A \backslash X = {\bf P}^1$ wildly branched above two 
points $P,Q$ and unramified elsewhere, and suppose that there exists a subgroup $G$ of $A$ such that 
$X \rightarrow G \backslash X = {\bf P}^1$ is branched above two points $x,y$, with $x$ tamely ramified and $y$ wildly ramified, such that $x$ and $y$ map to $P$ and $Q$ respectively on $A \backslash X$. Then $|A_{0,P}| = |A_{1,P}|(|A_{1,P}|-1)$.

The same holds if $\pi_A$ is wildly branched above one point $P$, branched of order 2 above two more points, and unramified elsewhere, if $x$ and $y$ map to $P$ on $A \backslash X$.}

\vv

\Proof Since $X$ is ordinary, the second ramification groups in both covers are
trivial. Let $g$ be the genus of $X$. The Riemann-Hurwitz formula applied to $\pi_G$ gives that 
$$ M:= \frac{2(g-1)}{|G|} = \frac{|G_{1,y}|-2}{|G_{0,y}|}-\frac{1}{|G_{0,x}|}, $$
so that $0<M<1$. Applying similarly Hurwitz's formula to $\pi_A$ and dividing $M$ by 
the result gives that 
$$[A:G] = M \frac{|A_{0,P}|\cdot |A_{0,Q}|}{|A_{0,Q}|(|A_{1,P}|-2)+|A_{0,P}|(|A_{1,Q}|-2)}. $$ 
Since both terms in the denominator are positive (as $P$ and $Q$ are wildly ramified and $p>2$),  
$$ [A:G] \leq M \frac{|A_{0,P}|}{|A_{1,P}|-2}. $$
Henn has proved (\cite{Henn:78}, Lemma 1) that $\frac{|A_{0,P}|}{|A_{1,P}|}$ divides $|A_{1,P}|-1$, so unless if $|A_{0,P}|=|A_{1,P}|(|A_{1,P}|-1)$, we have $[A:G] \leq M |A_{1,P}|$. Since ramification above $x$ in $\pi_G$ is only tame, the $p$-part $|A_{1,P}|$ is bounded by $[A:G]$, so that finally, $[A:G] \leq M [A:G]$, 
a contradiction to $M<1$. 

A similar computation works for the second case. \ee

\vv

(10.4) {\sl Proof of proposition 4.} \ \ We see from (10.1) that, except for $q=2,3$, $X({\frak n})$ has more than $12(g(X({\frak n}))-1)$ 
automorphisms. 

Let $\Gamma$ be the Schottky group of $X({\frak n})$ and $N$ its normalizer. We observe that $G({\frak n}) \subseteq A:=\A(X({\frak n}))$, and the cover $X({\frak n}) \rightarrow X(1)$ is a cover of ${\bf P}^1$ by a Mumford curve ramified above exactly two points $P,Q$ with ramification indices 
$q+1$ and $q^d(q-1)$ (by (10.1)). But such covers where classified in (4.6), and the only case in which the ramification behaviour is compatible with the one
of $X({\frak n}) \rightarrow X(1)$ is case $(A1)$ with $P=PGL,t_1=dt$. Since we have
a tower of coverings of the form $X({\frak n}) \rightarrow X(1) \rightarrow A \backslash X({\frak n})$, this implies that $N$ contains $N' = PGL(2,q) \ast_{B(t,n_-)} B(td,n_-)$. Again, the cover $\pi_A$ has to belong to one of the cases in (4.6) and (4.7), since $X({\frak n})$ has more than $12(g(X({\frak n}))-1)$ automorphisms (for $q>3$).
The only possible such $N$ containing $N'$ are $(A1)$ and $(F1)$ with $t_1=t_2=0$ (the latter since the bound $12(g(X({\frak n}))-1)$ is exceeded). By the previous lemma (since $p \neq 2$), we conclude that $|A_{0,P}| = |A_{1,P}|(|A_{1,P}|-1).$ In case $(A1)$, we have $|A_{0,P}|=p^{t_2}n_+$, and the above identity leads to $t=t_2$. On the other hand, $2t|t_2$; a contradiction. Similarly, in case $(F1)$ with $t_1=t_2=0$, we have $|A_{0,P}|=p^{t_3}n_-$, leading to $t=t_3$. On the other hand, $dt|t_3$ for $d>1$ since $N$ contains
$N'$; a contradiction. We conclude that $N=N'$ and $\A(X({\frak n}))=G({\frak n})$. 

Application of (7.4) to this situation now leads to the fact that closed the stratum 
of Mumford curves of genus $g(X({\frak n}))$ having an automorphism group containing $G({\frak n})$ (hence, $N$ as
normalizer of their Schottky group) has dimension zero. One of the points in this finite set corresponds to the Drinfeld modular curve above; this proves the final claim of proposition 4. \ee

\vv

{\bf (10.5) Remark.} \ \ The curve $X({\frak n})$ is defined over a finite extension $F_{\frak n}$ of $F$, generated by the ${\frak n}$-torsion of the Carlitz module over $F$. A standard argument in model
theory implies that for all but a finite number of primes $\frak p$ of $F_{\frak n}$, the curve $X({\frak n}) \times_{F_{\frak n}} F_{\frak n}/{\frak p}$ has the same automorphism group. It would be interesting to know what happens at such special primes ({\sl cf.} Adler \cite{Adler:97}, Rajan \cite{Rajan:99}). 

\vv

\vv

{\footnotesize

\noindent {\bf Acknowledgments.} The first author is post-doctoral fellow of the Fund for Scientific Research - Flanders (FWO -Vlaanderen). The third author is
supported by  the TMR-programme of the EU. We thank F.\ Herrlich for sending us a copy of \cite{Herrlich:78} and \cite{Herrlich:84}, and for a constructive discussion concerning (7.6). This work was done at the MPIM; the authors are delighted to notice that the emblem of this institution is an icosahedron that goes particularly well with proposition 2.}

\vv

\vv

\section*{Correction} 

The main theorem of \cite{CKK}, claiming to give an upper bound for the number of automorphisms of a Mumford curve in characteristic $p>0$, is false. The correct upper bound, due to van der Put and Voskuil, is the main theorem of \cite{vdPV}. The errors  in \cite{CKK} are explained in Remarks 4.8 and 8.27 of \cite{vdPV}. This also implies that the subsequent study, in \cite{CK1} and \cite{CK2}, of curves attaining the (false) upper bound from \cite{CKK} does not determine all Mumford curves with maximal automorphism group w.r.t.\ their genus; see \cite{vdPV} for corrections. 

The further main results from the introduction of \cite{CKK} remain valid in the following sense: 
\begin{enumerate}
\item Proposition 3 from \cite{CKK} describes the family of ``Artin-Schreier-Mumford curves" that indeed attains the upper bound from the main theorem of \cite{CKK}; but this bound is not the maximal one. 

\item Proposition 4 of \cite{CKK} describing the automorphism group of Drinfeld modular curves for principal congruence groups in characteristic $>3$ remains correct, even with the new classification from \cite{vdPV} (observing, as described in the note added in proof to \cite{CKK}, that the curves are not equivariantly deformation rigid, see \cite{CK}). 
\end{enumerate} 

\footnotesize{

\bibliographystyle{amsplain}

\vv

\noindent Max-Planck-Institut f{\"u}r Mathematik, P.O.\ Box 7280,
D-53072 Bonn (gc\footnote{correspondence should be sent to this author at this address}, fk \& ak)

\vv

\noindent Ghent University, Department of Pure Mathematics, 
Galglaan 2, B-9000 Ghent (gc)

\vv

\noindent Graduate School of Mathematics, Kyushu University 33, Fukuoka 812-8581, Japan (fk)

\vv

\noindent e-mail: {\tt gc@cage.rug.ac.be, fkato@math.kyushu-u.ac.jp, kontogeo@mpim-bonn.mpg.de}
}

\vv

\vv

{\bf Note added in proof.} \ \ The final sentence of proposition 4 
is incorrect. This is due to the fact that theorem (7.4), quoted
from [14], can fail to hold in positive characteristic. However, 
it does apply in the case of proposition 2 and 3 -- see forthcomming
work of the first two authors.

\end{document}